\documentclass{amsart}
\usepackage{mathpazo}
\usepackage{amssymb,color}
\usepackage{graphicx,epsfig}
\usepackage{amsmath}
\graphicspath{{figures/}}

\newtheorem{definition}{Definition}[section]
\newtheorem{thrm}{Theorem}[section]
\newtheorem{lmm}{Lemma}[section]
\newtheorem{crllr}{Corollary}[section]
\newtheorem{prpstn}{Proposition}[section]

\newtheorem{xmpl}{Example}[section]
\newtheorem{clm}{Claim}[section]

\theoremstyle{remark}
\newtheorem{rmrk}{Remark}[section]

\newcommand{\sgn}{{\rm sgn}\kern 0.12em}
\newcommand{\argmin}{{\rm argmin}\kern 0.12em}
\newcommand{\cC}{{\mathcal C}}
\newcommand{\cE}{{\mathcal E}}
\newcommand{\cL}{{\mathcal L}}

\newcommand{\cH}{{\mathcal H}}
\newcommand{\cT}{{\mathcal T}}
\newcommand{\cS}{\mathcal S}
\newcommand{\cO}{\mathcal O}
\newcommand{\cI}{\mathcal I}
\newcommand{\cU}{\mathcal U}

\newcommand{\eps}{\varepsilon}
\newcommand{\BB}{\mathbb{B}}
\newcommand{\R}{\mathbb{R}}
\newcommand{\N}{\mathbb{N}}
\newcommand{\bd}{{\rm bd}\kern 0.12em}
\newcommand{\inte}{{\rm int}\kern 0.12em}
\newcommand{\demi}{\frac{1}{2}}
\newcommand{\ie}{{\it i.e.}\,\,}
\newcommand{\cf}{{\it cf.}\,\,}

\newcommand{\barre}{\overline}
\begin{document}
\title[On second order differential equations with asymptotically small dissipation]{On the long time behavior
of second order differential equations with asymptotically small dissipation}

\author{Alexandre Cabot}
\address{D\'epartement de Math\'ematiques, Universit\'e Montpellier II, CC 051\\
Place Eug\`ene Bataillon, 34095 Montpellier Cedex 5, France}
\email{acabot@math.univ-montp2.fr}

\author{Hans Engler}
\address{Department of Mathematics, Georgetown University\\
Box 571233\\ Washington, DC 20057\\
USA}
\email{engler@georgetown.edu}

\author{S\'ebastien Gadat}
\address{Laboratoire de Statistique et Probabilit\'es, Universit\'e Paul Sabatier\\
31062 Toulouse Cedex 9, France}
\email{Sebastien.Gadat@math.ups-tlse.fr}

\subjclass{34G20, 34A12, 34D05}

\keywords{Differential equation, dissipative dynamical system,
 vanishing damping, averaged gradient system, asymptotic behavior, Bessel equation}

\begin{abstract}
We investigate the asymptotic properties as $t\to \infty$ of the following differential equation
in the Hilbert space $H$
$$\ddot{x}(t)+a(t)\dot{x}(t)+ \nabla G(x(t))=0, \quad t\geq 0, \leqno (\cS)$$
where the map $a:\R_+\to \R_+$ is non increasing and the potential $G:H\to \R$
 is of class $\cC^1$. If the coefficient  $a(t)$ is constant and positive, we recover
the so-called ``Heavy Ball with Friction'' system.  On the other hand, when $a(t)=1/(t+1)$ we obtain the trajectories associated to some averaged gradient system.
 Our analysis is mainly based on the existence
of some suitable energy function. When the function $G$ is convex, the condition
$\int_0^\infty a(t) \, dt =\infty$ guarantees that the energy function converges toward
its minimum.
The more stringent condition $\int_0^{\infty} e^{-\int_0^t a(s)\, ds}dt<\infty$ is necessary
 to obtain the convergence of the trajectories of $(\cS)$
toward some minimum point of $G$. In the one-dimensional setting,
 a precise description of the convergence of solutions is given for a general non-convex
 function $G$. We show that in this case the set of initial conditions for which solutions converge to a local minimum is open and dense.
\end{abstract}

\maketitle
\section{Introduction}
Throughout this paper, we study the differential equation
\begin{equation*} \ddot{x}(t)+a(t)\dot{x}(t)+ \nabla G(x(t))=0, \quad t\geq 0 \leqno (\cS)
\end{equation*}
in a finite- or infinite-dimensional Hilbert space $H$, where the map $G:H \to \R$ is at least of
class $\cC^1$ and $a:\R_+\to \R_+$ is a non increasing function. To motivate our study, let us describe four examples and applications
which are intimately connected with equation $(\cS)$. \\

{\em Averaged gradient system} \,For the potential $G$, the much studied
gradient flow is defined as the solution map $y(0) \mapsto y(s),\, s \ge 0$ of
the differential equation
\[\dot{y}(s) =- g(y(s)) = -\nabla G(y(s)) \, .\]
It is of interest to consider the case where $\dot{y}(s)$ is proportional, not
to the instantaneous value of $\nabla G(y(s))$, but to some average of $\nabla
G(y(\tau)), \, \tau \le s$. The simplest such equation is
\begin{equation} \label{average}
\dot{z}(s) + \frac{1}{s} \int_0^s g(z(\tau)) d\tau = 0.
\end{equation}
For more general gradient systems with memory terms involving kernels, we refer for example
to \cite{GouMun}.
After multiplying equation (\ref{average}) by $s$ and differentiating, this leads to the
ordinary differential equation
\begin{equation*} s \ddot{z}(s) + \dot{z}(s) + g(z(s)) = 0
\end{equation*}
which becomes
\begin{equation} \label{bessel0} \ddot{x}(t) + \frac{1}{t}\dot{x}(t) +
g(x(t)) = 0
\end{equation}
after the change of variables $s = \frac{t^2}{4},\, t = 2 \sqrt{s}, \,  x(t) =
z\left(\frac{t^2}{4}\right), z(s) = x(2 \sqrt{s})$. This is the problem $(\cS)$
with $a(t) = \frac{1}{t}$. We note that in the special case where $g(\xi) = \xi$, (\ref{bessel0}) is a Bessel equation. All
solutions with a finite limit at $t=0$ are multiples of
\[J_0(t) = \sum_{k=0}^\infty (-1)^k \frac{t^{2k}}{2^{2k}
(k!)^2} \,.\] It is well known that
\[J_0(t) \sim \sqrt{\frac{2}{\pi t}} \cos \left(t - \frac{\pi}{4}\right) \]
and therefore
\[z(s) = C\cdot J_0(2\sqrt{s}) \sim C \cdot \sqrt[4]{\frac{1}{ s }} \cos \left(2 \sqrt{s} - \frac{\pi}{4}\right) \]
for some suitable constant $C$ as $s, \,t \to \infty$.
Thus the solution $z$ of the averaged system
 (\ref{average}) converges to zero just as the solution $y(s) =y(0)\,e^{-s}$
 of the corresponding gradient system does, but it does so much
 more slowly (at an algebraic rate), and it oscillates infinitely often.
 Our work will generalize this simple famous example using several cases for
 $a$ and $G$. The case where $g(x) = x^3 - x$ and $H = \R$ was discussed in \cite{miles82}. \\

{\em Heavy Ball with Friction system}\,
A particular attention has been recently devoted to the so-called ``Heavy Ball with Friction'' system
$$\ddot{x}(t)+\gamma\, \dot{x}(t)+ \nabla G(x(t))=0,  \leqno (HBF)$$
where $\gamma>0$ is a positive damping parameter. From a mechanical point of view, the $(HBF)$ system
corresponds to the equation describing the motion of a material point subjected to the conservative force
$-\nabla G(x)$ and the viscous friction force $-\gamma\, \dot{x}$. \\
The $(HBF)$ system is dissipative
and can be studied in the classical framework of the theory of dissipative dynamical systems (\cf
Hale \cite{Hal}, Haraux \cite{Har}).
The presence of the inertial term $\ddot x(t)$  allows to overcome some
drawbacks of the steepest descent method. The main interest of the $(HBF)$ system in numerical
optimization is that it is not a descent method: it permits to go up and down along the graph of $G$.
The trajectories of $(HBF)$ are known to be convergent toward
a critical point of $G$ under various assumptions like convexity, analyticity,~... In the convex
setting, the proof of convergence relies on the Opial lemma, see Alvarez \cite{Alv}, Attouch-Goudou-Redont
\cite{AttGouRed}, while it uses the
Lojasiewicz inequality in the case of analytic assumptions, (\cf Haraux-Jendoubi \cite{HarJen}).\\
In the above $(HBF)$ model, the damping coefficient $\gamma$ is constant. A natural extension consists
in introducing a time-dependent damping coefficient, thus leading to the system $(\cS)$. In our paper, we
will focus on the important case corresponding to a vanishing damping term $a(t)$, \ie $a(t)\to 0$
as $t\to \infty$. It is clear that the decay properties of the map $a$  play a central role
in the asymptotic behavior of $(\cS)$. In particular, if the quantity $a(t)$ tends to $0$ too rapidly as
$t\to \infty$, convergence of the trajectory may fail (think about the extreme case of $a\equiv0$ for instance).\\

{\em Semilinear Elliptic equations}\,
Consider the semilinear elliptic system
\[\Delta u(y) + g(u(y)) = 0\]
in $\R^m$, where $u:\R^m \to \R^n$ is the unknown function. Radial
solutions $u(y) = x(|y|)$ of this system lead to the ordinary
differential equation
\[\ddot{x}(r) + \frac{m-1}{r}\dot{x}(r) + g(x(r)) = 0\,.\]
There has been a large amount of work on this problem; see e.g. \cite{Ni04} for a recent overview.

{\em Stochastic Approximation algorithms}\, The classical stochastic algorithm
introduced by \cite{robbins-monro} is used in many fields of approximation
theory. This method is frequently used to approximate, with a random version of
the explicit Euler scheme, the behavior of the ordinary differential equation
$\dot{x}(t)=-g   (x(t))$. If we denote $(X^n)_{n \in \N}$ the random
approximations, $(\omega^n)_{n \geq 1}$ and $(\eta^n)_{n \geq 1}$ two auxiliary
stochastic processes, the recursive approximation is generally written as
\begin{equation}\label{sa}
\left\{
\begin{array}{ll}
        X^0 \in \R^d \\
        X^{n+1} = X^n-\eps_{n+1} g(X^n,\omega^{n+1}) + \eps_{n+1}
        \eta^{n+1},\qquad \forall t \in \N
    \end{array}
\right.
\end{equation}
where the gain of the algorithm $\eps_{n}$ is a sequence of positive real
numbers and $\eta^{n}$ is a small residual perturbation which is zero in many
cases. Defining by $(\mathcal{F}_n)_{n \geq 1}$ the set of measurable events at
time $t=n$, solutions of (\ref{sa}) are shown to asymptotically behave like
those of the determinist o.d.e. $\dot{x}(t)=-g(x(t))$ provided $\eta^n =
o(\eps^n)$, $\Delta M^n = g(X^n) - g(X^n,\omega^{n+1})$ is the increment of a
(local) $\mathcal{F}_n$-martingale, and the sequence $(\eps_n)_{n \geq 1}$
satisfies the baseline assumptions:
$$
\sum_{n=1}^{\infty} \eps_n = \infty  \qquad \text{and} \qquad
\sum_{n=1}^{\infty} \eps_n^{1+\alpha} < \infty, \quad \mbox{ for some }\alpha
>0.
$$
The assumption on the martingale increment implies that $ g(X^n) = \mathbb E
\left[ g(X^n,\omega^{n+1}) | \mathcal{F}^n \right]. $ A very common case occurs
when $(\omega^n)_{n \ge 1}$ is a sequence of independent identically
distributed
 variables with distribution $\mu$ and $g(x,.)$ is $\mu$-integrable:
$$
\forall x \qquad  g(x):=\int g(x,\omega) \mu(d\omega).
$$

 This yields the stochastic gradient descent algorithm when $g$ is the gradient operator
of a potential $G$. One recent application \cite{gadat-younes} applied this stochastic
gradient to perform feature selection among a large amount of variables with a
simpler Euler scheme
$
X^{n+1} = X^n - \eps_n g(X^n,\omega^{n+1}).
$
 Further developments have shown that in some cases, the random variable
 $g(X^n,.)$ may have a large variance and the stochastic approximation of
$g(X^n)$ by $g(X^n,\omega^{n+1})$ can be numerically improved using the
following modified recursive definition:
\begin{equation}\label{sa-modified}
        X^{n+1} = X^n-\eps_{n+1} \frac{\displaystyle\sum_{i=1}^n
        \eps_i g(X^i,\omega^{i+1})}{\displaystyle\sum_{i=1}^n \eps_i}.
\end{equation}
One can think about (\ref{sa-modified}) as a way to improve the instability of
the gradient estimate $g(X^n)$ by an average on the variables
$\{g(X^k,\omega^{k+1}), k \leq n\}$ whose weights correspond to the $\eps_n$.
Actually, one can show (see the proof in  Appendix A) that the limit o.d.e. is
 given by an equation of type $(\cS)$:
\begin{equation}
\ddot{X}(t)= -
\frac{\dot{X}(t)+g(X(t))}{t+\beta},
\end{equation}
for some $\beta\geq 0$. In the particular case $\beta=0$, we obtain
the average gradient system equation (\ref{average}).\\

The analysis of the asymptotic behavior of $(\cS)$ is based on the use of the energy function
$\cE$ defined by $\cE(t) =\frac{1}{2}|\dot{x}(t)|^2 + G(x(t)) $ for every $t\geq 0$.
Under convex-like assumptions on $G$, we prove the convergence of the
quantity $\cE (t)$ toward $\min G$ as $t\to \infty$, provided that $\int_0^\infty a(t)\,
dt =\infty$. This condition expresses that the damping coefficient $a(t)$ slowly tends to $0$
as $t\to \infty$. Such a condition has already been pointed out for the steepest descent
method combined with Tikhonov viscosity-regularization in convex minimization \cite{AttCom96}
 as well as for the stabilization of nonlinear oscillators \cite{AttCza,Cab}.
When the convex function $G$ has a unique minimum $\barre{x}$, condition $\int_0^\infty a(t)\,
dt =\infty$ is sufficient to ensure the convergence of the trajectories of $(\cS)$ toward
$\barre{x}$. If the function $G$ has a set of non-isolated equilibria, the more stringent condition $\int_0^{\infty} e^{-\int_0^t a(s)\, ds}dt<\infty$ is necessary
 to obtain the convergence of the trajectories of $(\cS)$
toward some minimum point of $G$. Notice that the previous condition fails if
$a(t)=\frac{1}{t+1}$ for every $t\geq 0$, which shows that the averaged gradient system
defined above is divergent when the convex function $G$ has multiple minima.\\
We also have substantial results in the non-convex setting, when the function $G$
 has finitely many critical points. Under the slow condition $\int_0^\infty a(t)\,dt =\infty$, we then prove that
 the energy function $\cE (t)$ converges toward a critical value as $t\to \infty$. If moreover
there exists $c>0$ such that $a(t)\geq \frac{c}{t+1}$ for every $t\geq 0$, we show that a
Cesaro average of the solution $x$ converges toward some critical point of $G$. Finally,
in the one-dimensional setting,
 a precise description of the convergence of solutions is given for a general non-convex
 function $G$. We show that in this case the set of initial conditions for which solutions
converge to a local minimum is open and dense.\\

{\em Outline of the paper }\,
Our work starts with a global existence result of solutions to $(\cS)$, based on the use of the Lyapounov function $\cE$. Section \ref{se.energy_estim} is concerned with the asymptotic
behavior of the energy function $\mathcal{E}$ under
convex-like hypotheses on $G$, and provides estimates on the speed of convergence
of the quantity $\mathcal{E}(t)$ toward $\inf G$ as $t\to \infty$.
 Section \ref{ConvTraj} explores the convergence of the trajectories
of $(\cS)$ in the general setting of convex functions having multiple minima.
In section \ref{sec:NonConv}, we study the asymptotic behavior of $(\cS)$
 in the non-convex case when $G$ has finitely many critical points.
Finally section \ref{sec:OneDim} is dedicated to the very special one
dimensional case. Details for the stochastic gradient descent algorithm are given in appendix A, and some special equations are discussed in appendix B.

\section{General Facts}
In the entire paper, we will denote
by $G$ a $\cC^1$ potential map from an Hilbert space $H$ into $\R$ for which
the gradient $g = \nabla G$ is Lipschitz continuous, uniformly on bounded sets.
 Given a function $a:\R_+ \to \R_+$, we will consider
the following dynamical system
\begin{equation*} \ddot{x}(t)+a(t)\dot{x}(t)+ g(x(t))=0, \quad t\geq 0 \,. \leqno (\cS)
\end{equation*}
Let us start with a basic result on existence and uniqueness for solutions of $(\cS)$.
In the next statement the map $a$ may have a
singularity at $t=0$ so as to cover cases like $a(t)=1/t$, for $t>0$.

\begin{prpstn}\label{pr.existence}
(a) Suppose $a:(0,\infty) \to \R_+$ is continuous on $(0,\infty)$
and integrable on $(0,1)$. Then for any $(x_0, \, x_1) \in H \times H$, there exists a unique solution $x(\cdot) \in \cC^2([0,T),H)$ of $(\cS)$ satisfying $x(0) = x_0,\, \dot{x}(0) = x_1$ on some maximal time interval $[0,T) \subset [0,\infty)$. \\
(b) Suppose $a:(0,\infty) \to \R_+$ is continuous and there exists $c>0$ such that
$a(t) \le \frac{c}{t}$ for $t\in (0,1]$. Then for any $x_0 \in H$, there exists a unique solution $x(\cdot) \in \cC^2((0,T),H)\cap \cC^1([0,T),H)$ of $(\cS)$ satisfying $x(0) = x_0,\, \dot{x}(0) = 0$ on some maximal time interval $[0,T) \subset [0,\infty)$.
\end{prpstn}
The previous proposition can be proved with standard arguments for ordinary differential equations.
The result below states the decay property of the energy function $\cE$ defined by
\begin{equation}\label{eq.define_cE}
\cE(t) =\frac{1}{2}|\dot{x}(t)|^2 + G(x(t)) \,.
\end{equation}
A global existence result is then derived when the potential function $G$ is bounded from below.
The existence of the Lyapounov function $\cE$ will be a crucial tool for the analysis of the
asymptotic behavior of $(\cS)$.
\begin{prpstn}\label{pr.general}
Let $a: \R_+ \to \R_+$ a continuous map and  let $G: H\to \R$ a function of class $\cC^1$ such
that $\nabla G$ is Lipschitz continuous on the bounded sets of $H$.
Let $x$ be a solution to $(\cS)$ defined on some interval $[0,T)$, with $T\leq \infty$.
\\ (a) For every $t\in [0,T)$, the following equality holds
\begin{equation} \label{energy1}
\frac{d}{dt} \cE(t) = -a(t)|\dot{x}(t)|^2
\end{equation}
and therefore for $0 \le s < t <T$
\begin{equation}
\label{energy2} \cE(s)- \int_s^t a(\tau) |\dot{x}(\tau)|^2 d \tau = \cE(t) \, .
\end{equation}
\\(b) If in addition $G$ is bounded from below on $H$, then
\begin{equation}
\label{energy3} \int_0^T a(t)|\dot{x}(t)|^2 dt < \infty \,
\end{equation}
and the solution exists for all $T > 0$.
\\(c) If also $G$ is coercive\footnote{Let us recall that the coercivity of $G$ means that
$G(\xi) \to \infty$ as $|\xi| \to \infty$.}, then all solutions to $(\cS)$ remain bounded
together with their first and second derivatives for all $t > 0$. The bound
depends only on the initial data.
\end{prpstn}

\begin{proof} (a) Equation (\ref{energy1}) follows by taking the scalar product of
 $(\cS)$ against $\dot{x}(t)$, and (\ref{energy2}) follows by integrating.

\noindent (b) If $G$ is bounded from below, then equality (\ref{energy1}) shows
that $t\mapsto \mathcal{E}(t)$ is decreasing and remains bounded.
 Estimate (\ref{energy3}) is then a consequence of equality (\ref{energy2}), and it also
follows that $\sup_{t<T} \cE(t) < \infty$. Therefore $\dot{x}$ is uniformly
bounded on $[0,T)$. If $T < \infty$, then the solution $x$ together with its
derivative has a limit at $t=T$ and therefore can be continued. Thus the
solution $x(t)$ exists for all $t$ and $\dot{x}$ is uniformly bounded by
quantities depending on the initial data.

\noindent (c) Using the coercivity of $G$ and the inequality $G(x(t)) \le \cE(t)\leq \cE(0)$,
 we derive that the map $x$ is uniformly bounded. Then also $\ddot{x}(t)$ is uniformly bounded due to the differential equation $(\cS)$ and this bound depends only on the initial data.
\end{proof}

If $a$ does not decrease to $0$ too rapidly, then the derivative of any solution of ($\cS$)
must be arbitrarily small on arbitrarily long time intervals, infinitely often.

\medskip
\begin{prpstn}\label{pr.slowdown}
Let $a:\R_+ \to \R_+$ be a non increasing map such that $\int_0^\infty a(s) ds = \infty$. Let
$G:H\to \R$ be a coercive function of class $\cC^1$ such that $\nabla G$ is Lipschitz
 continuous on the bounded sets of $H$. Then, any solution $x$ to the differential
equation $(\cS)$ satisfies, for every $T>0$,
\[\liminf_{t \to \infty} \sup_{s \in [t,t+T]} |\dot{x}(s)| = 0 \, .\]
\end{prpstn}

\begin{proof} Suppose not, then there exist $\eps > 0, \, T>0$
such that for all $k$
\[\sup_{s \in [kT, (k+1)T]}|\dot{x}(s)| > \eps \]
and thus there are $t_k \in [kT, (k+1)T]$ such that
$|\dot{x}(t_k)|>\eps$. Since the map $\dot x$ is Lipschitz continuous, there exists
some fixed $\delta>0$ such that $|\dot{x}(t)| \ge \eps/2$ on
$[t_k-\delta,t_k+\delta]$ for every $k\in \N$. Since the map $a$ is non increasing, we have
\begin{eqnarray*}
\sum_k a((k+1)T) \le \sum_k a(t_k) \le \sum_k \frac{4}{\delta
\eps^2} \int_{t_k-\delta}^{t_k} a(t)|\dot{x}(t)|^2dt\le \frac{4}{\delta \eps^2} \int_0^\infty a(t)|\dot{x}(t)|^2dt.
\end{eqnarray*}
Recalling that $\int_0^\infty a(t)|\dot{x}(t)|^2dt \leq \cE(0) -\inf G<\infty$,
we infer that $\int_0^\infty a(t) dt < \infty$, a contradiction.
\end{proof}

We next show that if $\dot{x}$ is small on some interval, then $g(x(t))$
is proportionally small on a slightly shorter interval. This implies that if solutions
slow down for a long time interval, they must be near a critical point of $G$.

\begin{prpstn}\label{pr.g.estimate}
Let $a:\R_+ \to \R_+$ be a non increasing map. Let
$G:H\to \R$ be a function of class $\cC^1$ such that $g=\nabla G$ is Lipschitz continuous
on the bounded sets of $H$. If $x$ is a solution of ($\cS$)
and $|\dot{x}(t)| \le \eps$ on $[T_0,T_1]$, then for every $\delta \in \left]0,
\frac{T_1-T_0}{2}\right]$, we have
\[\forall t\in [T_0+\delta, T_1-\delta], \qquad |g(x(t))| \le \left(\frac{2}{\delta} + a(T_0) + \frac{L\delta}{2}\right) \eps,\]
where $L>0$ is a Lipschitz constant of the map $g$ on the set $x([T_0,T_1])$.
\end{prpstn}

\begin{proof} Suppose $|\dot{x}(t)| \le \eps$ on $[T_0,T_1]$ with
$T_1 - T_0 \ge 2\delta$. Let the map $t\mapsto |g(x(t))|$ attain its maximum
on $[T_0+\delta,T_1-\delta]$ at $t=t_0$. Since the map $g$ is $L$-Lipschitz continuous
on the set $x([T_0,T_1])$,  we  have for every $t\in [t_0-\delta, t_0+\delta]$
\[|g(x(t_0)) - g(x(t))| \le L|x(t_0)-x(t)| \le L \eps|t_0-t|\]
 and thus
\[|\ddot{x}(t) + g(x(t_0))| = |a(t)\dot x(t) + g(x(t)) - g(x(t_0))| \le a(T_0)\eps + L \eps
|t_0-t|\]
for the same range of $t$. We then have
\begin{eqnarray*}
\eps \ge |\dot x(t)| &\ge& \left|\int_{t_0}^t \ddot{x}(s)ds \right| - |\dot x(t_0)| \\
&\ge& \left|\int_{t_0}^t g(x(t_0)) ds\right| - \int_{t_0}^t|\ddot{x}(s) + g(x(t_0))| ds
-\eps \\
&\ge& |t_0-t||g(x(t_0))| - \eps - a(T_0) \eps|t_0-t| - \frac{L}{2}
\eps|t_0-t|^2
\end{eqnarray*}
and therefore
\[|g(x(t_0))| \le \left(\frac{2}{|t_0-t|} + a(T_0) + \frac{L}{2}|t_0-t| \right)
\eps  \,.\] Set $t = t_0 \pm \delta$ to conclude $|g(x(t_0))| \le
\left(\frac{2}{\delta} + a(T_0)+ \frac{L \delta}{2} \right) \eps$.
\end{proof}

By combining Propositions \ref{pr.slowdown} and \ref{pr.g.estimate}, we derive the
following corollary.
\begin{crllr}
Under the assumptions of Proposition \ref{pr.slowdown},  any solution $x$ to the differential
equation $(\cS)$ satisfies, for every $T>0$,
\[\liminf_{t \to \infty} \sup_{s \in [t,t+T]} |g(x(s))| = 0 \, .\]
\end{crllr}
\noindent The proof is immediate and left to the reader.\\

The last result of this section establishes that,
 if a solution $x$ to $(\cS)$ converges toward some $\barre{x}\in H$, then $\barre{x}$ is a
stationary point of $\nabla G$ and moreover the velocity $\dot x(t)$ and the acceleration
 $\ddot x(t)$ tend to $0$ as $t\to \infty$.
\begin{prpstn}\label{pr.stationary}
Let $a: \R_+ \to \R_+$ be a bounded continuous map and let $G: H\to \R$
be a function of class $\cC^1$  such
that $\nabla G$ is Lipschitz continuous on the bounded sets of $H$.
 Consider a solution $x$ to $(\cS)$ and assume that
there exists $\barre{x}\in H$ such that $\lim_{t\to \infty}x(t)=\barre{x}$. Then we have
$\lim_{t\to \infty}\dot x(t)=\lim_{t\to \infty}\ddot x(t)=0$ and the vector $\barre{x}$ satisfies
$\nabla G(\barre{x})=0$.
\end{prpstn}
\begin{proof} Since $x(\cdot)$ converges, it is uniformly bounded. Due to the inequality
$\cE(t)\leq \cE(0)$ for every $t\geq 0$, $\dot{x}(\cdot)$ is also uniformly bounded, and just as in the proof of Proposition \ref{pr.general} (c), the map $\ddot x(\cdot)$ is uniformly bounded as well by some constant $M>0$.
Landau's inequality applied to the map $t\mapsto x(t)-\barre{x}$ yields, for every $t\geq 0$
$$\sup_{[t,\infty[}|\dot x|\leq 2\sqrt{\sup_{[t,\infty[}|x-\barre{x}|\,.\,
\sup_{[t,\infty[}|\ddot x|}\leq
2\sqrt{M}\,\sqrt{\sup_{[t,\infty[}|x-\barre{x}|}.$$ By using the assumption
$\lim_{t\to \infty}x(t)=\barre{x}$ and letting $t\to \infty$ in the above
inequality, we derive that $\lim_{t\to \infty}\dot x(t)=0$. Since $\lim_{t\to
\infty}\nabla G(x(t))=\nabla G(\barre{x})$, the differential equation $(\cS)$
shows that $\lim_{t\to \infty}\ddot x(t)=-\nabla G(\barre{x})$. If $\nabla
G(\barre{x})\neq~0$, an immediate integration gives the equivalence $\dot
x(t)\sim -t\,\nabla G(\barre{x})$ as $t\to \infty$, a contradiction. Thus we
conclude that $\nabla G(\barre{x})=0$ and $\lim_{t\to \infty}\ddot x(t)=0$.
\end{proof}

\section{Case of a convex-like potential. Energy estimates}\label{se.energy_estim}
As in the previous section, $(\cS)$ is studied on a general Hilbert space $H$ in this section.
Throughout, we will assume that the function $G$ satisfies the following condition:
 there exist $z\in\argmin G$ and $\theta\in \R_+$ such that
\begin{equation}\label{eq.ineg_base}
\forall x\in H, \qquad G(x)-G(z)\leq \theta\, \langle \nabla G(x), x-z\rangle.
\end{equation}

This can be viewed as a generalization of the notion of convexity, and it also generalizes Euler's identity
for homogeneous functions. Indeed, if
$G$ is convex,
 then condition (\ref{eq.ineg_base}) is satisfied with $\theta=1$ for every $z\in H$.
Now assume that $G$ is defined by $G(x)=\demi\, \langle Ax,x\rangle$ where
$A\in \cL(H)$ is symmetric and positive. We then have $g(x)=Ax$ and inequality
(\ref{eq.ineg_base}) holds as an equality with $\theta=\demi$ and $z=0$.
Finally, if $G$ is defined by $G(x)=\frac{|x|^p}{p}$ with $p>1$, we have $g(x)=
x\, |x|^{p-2}$ and inequality (\ref{eq.ineg_base}) is satisfied with
$\theta=\frac{1}{p}$ and $z=0$ as an equality.

\subsection{A result of summability}
First we give a result of summability of the function $t\mapsto \cE(t)-\min G$
over $\R$, with respect to some measure depending on the map $a$. This property
 will imply some convergence results
 on $\cE$  provided some weak hypotheses on the function $a$.
\begin{prpstn} \label{th.ener_L1}
Let $a: \R_+ \to \R_+$ be a non increasing and differentiable map.
 Let $G: H\to \R$ be a coercive function of class $\cC^1$  such
that $\nabla G$ is Lipschitz continuous on the bounded sets of $H$.
Assume that $\argmin G\neq \emptyset$ and that there exist
 $z\in\argmin G$ and $\theta\in \R_+$ such that condition (\ref{eq.ineg_base}) holds.
 Then, any solution $x$ to the differential equation $(\cS)$ satisfies the following estimate
$$\int_0^{\infty}a(t) \, (\cE(t)-\min G)\, dt <\infty.$$
\end{prpstn}
\begin{proof}
Let us define the function $h:\R_+\to \R$ by
$$h(t)=\frac{a(t)}{2}\, |x(t)-z|^2+\langle \dot{x}(t), x(t)-z\rangle.$$
By differentiating, we find:
$$\dot{h}(t)= \frac{\dot{a}(t)}{2}\, |x(t)-z|^2+a(t)\,\langle \dot{x}(t), x(t)-z\rangle
+\langle \ddot{x}(t), x(t)-z\rangle+ |\dot{x}(t)|^2.$$ Since
$\dot{a}(t) \leq 0$, we derive that
\begin{eqnarray*}
\dot{h}(t)&\leq&|\dot{x}(t)|^2+\, \langle
\ddot{x}(t)+a(t)\,\dot{x}(t), x(t)-z\rangle
\nonumber\\
&\leq&|\dot{x}(t)|^2-\, \langle g(x(t)), x(t)-z\rangle.
\end{eqnarray*}
Let us now fix $m\in \left]0,\frac{1}{\theta+1/2}\right]$. Recalling that
$\dot{\cE}(t)=-a(t)\, |\dot{x}(t)|^2$, we find
\begin{small}
$$\hspace{-6cm} \dot{\cE}(t)+m\, a(t)\, (\cE(t)-\min G)+\theta\, m\, a(t) \,\dot{h}(t)  = $$
$$\hspace{1cm}(-1+(\theta+1/2)\,m)\,a(t)\, |\dot{x}(t)|^2  +m\, a(t)\, \left[G(x(t))-\min G -\theta\,
  \langle g(x(t)),x(t)-z\rangle\right].$$
\end{small}
Using condition  (\ref{eq.ineg_base}) and the fact that $m\leq \frac{1}{\theta+1/2}$,
we deduce
 \begin{equation}\label{eq.relation_E_dotE}
\dot{\cE}(t)+m\, a(t)\, (\cE(t)-\min G)+\theta\, m\, a(t)
\,\dot{h}(t)\leq 0.
\end{equation}
Let us integrate the previous inequality on $[0,t]$. Since
$\cE(t)\geq \min G$, we obtain
\begin{equation}\label{eq.majo_energie}
m\,\int_0^{t}a(s) \, (\cE(s)-\min G)\, ds\leq \cE(0)-\min
G\,-\theta\, m\,
 \int_0^t a(s)\dot{h}(s)\, ds.
\end{equation}
Then, remark that
\begin{equation}\label{eq.integ_parts}
\int_0^t a(s)\dot{h}(s)\, ds=a(t)h(t)-a(0)h(0)-\int_0^t
\dot{a}(s)\,{h}(s)\, ds.
\end{equation}
From the decay of the energy function $\cE$, it ensues that $t\mapsto
|\dot{x}(t)|$ and $t\mapsto G(x(t))$ are bounded. Since the map $G$ is
coercive, we infer that the map $t\mapsto |x(t)|$ is bounded. From the
expression of $h$, and the boundedness of $t\mapsto \mathcal{E}(t)$ and thus of
$t\mapsto |\dot{x}(t)|$, we immediately conclude the existence of $M>0$ such
that $|h(t)|\leq M$ for every $t\geq 0$. We then derive from
 (\ref{eq.integ_parts}) that
\begin{eqnarray*}
\left| \int_0^t a(s)\dot{h}(s)\, ds\right|&\leq& M\, a(t)+M\, a(0)+
M\, \int_0^t |\dot{a}(s)|\, ds\\
&=&M\, a(t)+M\, a(0)+ M(a(0)-a(t))=2\, M\, a(0).
\end{eqnarray*}
From (\ref{eq.majo_energie}), we now have that
$$\forall t\geq 0, \qquad m\int_0^{t}a(s) \, (\cE(s)-\min G)\, ds\leq
\cE(0)-\min G+2\,\theta\, m\,M\, a(0)$$ and we conclude that
$\int_0^{\infty}a(s) \, (\cE(s)-\min G)\, ds <\infty.$
\end{proof}

Now, we can prove the convergence of $\cE(t)$ toward $\min G$ as $t\to \infty$, provided that
 $\int_0^{\infty}a(t) \, dt=\infty.$ Notice that this assumption amounts to saying that the
quantity $a(t)$ slowly tends to $0$ as $t\to\infty$.

\begin{crllr}\label{co.energy_conv}
Under the hypotheses of Proposition \ref{th.ener_L1}, assume moreover that
$\int_0^{\infty}a(t) \, dt =\infty$. Then $\lim_{t\to
\infty} \cE(t)=\min G$. As a consequence, $\lim_{t\to
\infty} |\dot{x}(t)|=0$ and $\lim_{t\to \infty} G(x(t))=\min G$.
\end{crllr}
\begin{proof} Let us argue by contradiction and assume that $\lim_{t\to \infty} \cE(t)>\min G$.
This implies the existence of $\eta>0$ such that $\cE(t) -\min G
\geq \eta$ for every $t\geq 0$. We deduce that
$$\int_0^{\infty}a(t) \, (\cE(t)-\min G)\, dt\geq
 \eta \, \int_0^{\infty}a(t) \, dt=\infty.$$
 This yields a contradiction and we obtain the conclusions that $\lim_{t\to
\infty} |\dot{x}(t)|=0$ and $\lim_{t\to \infty} G(x(t))=\min G$.
\end{proof}

The next corollary precises the speed of convergence of $\cE$ toward $\min G$
under some assumption on the  decay of $t \mapsto a(t)$.

\begin{crllr} Under the hypotheses of Proposition \ref{th.ener_L1}, assume moreover that
there exists $m>0$ such that $a(t)\geq m/t$ for $t$ large enough. Then
$$\cE(t)-\min G =o\left( \frac{1}{t\,a(t)}\right)\qquad \mbox{as } t \to \infty.$$
\end{crllr}
\begin{proof} Since the functions $a$ and $\cE$ are non increasing and positive, it is immediate
that the map $t \mapsto a(t)\, (\cE(t)-\min G)$ is also non increasing. In
particular, we obtain
$$\int_{t/2}^t a(s)\, (\cE(s)-\min G)\, ds \geq \frac{t}{2}\, a(t)\, (\cE(t)-\min G).$$
Since $\int_0^{\infty}a(s) \, (\cE(s)-\min G)\, ds <\infty$, the left member
of the above inequality tends to $0$ as $t\to \infty$, which implies that
$\displaystyle\lim_{t\to \infty} t\,a(t)\,(\cE(t)-\min G)=0$.
\end{proof}

\subsection{Case of a unique minimum}
In view of the previous results, we are able to investigate the question of the
convergence of the trajectories in the case of a unique minimum. Studies with
several minima are more complicated and will be detailed in section
\ref{ConvTraj} (convex setting), section \ref{sec:NonConv} (non-convex setting)
and section \ref{sec:OneDim} (one-dimensional case).

\begin{prpstn} Let $a: \R_+ \to \R_+$ be a non increasing and differentiable map such that
$\int_0^{\infty}a(t) \, dt =\infty$.
Consider a map $\alpha:\R_+\to \R_+$ such that $\alpha(t_n) \to 0$
$\Longrightarrow$ $t_n \to 0$ for every sequence $(t_n)\subset
\R_+$. Let $G: H\to \R$ be a coercive function of class $\cC^1$ such
that $\nabla G$ is Lipschitz continuous on the bounded sets of $H$.
Given $\barre{x}\in H$, assume that
\begin{equation}\label{eq.mino_alpha}
\forall x\in H, \qquad G(x)\geq G(\barre{x})+\alpha(|x-\barre{x}|),
\end{equation}
and that there exists $\theta\in \R_+$ such that condition (\ref{eq.ineg_base}) holds with
$z=\barre{x}$. Then, any solution $x$ to the differential equation $(\cS)$ satisfies
 $\displaystyle\lim_{t\to \infty}x(t)=\barre{x}$ strongly in~$H$.
\end{prpstn}
\begin{proof} By applying Corollary \ref{co.energy_conv}, we obtain
$\lim_{t\to \infty} G(x(t))=\min G=G(\barre{x})$. From assumption (\ref{eq.mino_alpha}),
 we deduce that $\lim_{t\to \infty}\alpha(|x(t)-\barre{x}|)=0$ and we finally conclude
 that $ \lim_{t\to \infty} |x(t)-\barre{x}|=0$.
\end{proof}
If the stringent condition (\ref{eq.mino_alpha}) is not satisfied, one
can nevertheless obtain a result of weak convergence, as shown by the following statement.
\begin{prpstn}
 Let $a: \R_+ \to \R_+$ be a non increasing and differentiable map such that
$\int_0^{\infty}a(t) \, dt =\infty$.
 Let $G: H\to \R$ be a convex coercive function of class $\cC^1$
such that $\nabla G$ is Lipschitz continuous on the bounded sets of $H$.
If $\argmin G=\{\barre{x}\}$ for some $\barre{x}\in H$, then any solution
$x$ to the differential equation $(\cS)$ weakly converges to $\barre{x}$ in $H$.
\end{prpstn}
\begin{proof} Since $G$ is coercive, the trajectory $x$ is bounded. Hence there exist
$x_\infty \in H$ and a subsequence $(t_n)$ tending to $\infty$ such that
$\lim_{n\to \infty}x(t_n)=x_\infty$ weakly in $H$. Since $G$ is convex and continuous for the
strong topology, it is lower semicontinuous for the weak topology. Hence, we have:
$$G(x_\infty)\leq \liminf_{n\to \infty} G(x(t_n)).$$
On the other hand, by applying Corollary \ref{co.energy_conv}, we obtain
$ \lim_{t\to \infty} G(x(t))=\min G$. Therefore we deduce that
$G(x_\infty)\leq \min G$, \ie $x_\infty \in \argmin G=\{\barre{x}\}$. Hence
$\barre{x}$ is the unique limit point of the map $t\mapsto x(t)$ as $t\to
\infty$ for the weak topology. It ensues that $\lim_{t\to
\infty}x(t)=\barre{x}$ weakly in $H$.
\end{proof}

\subsection{Convergence rate of the energy function $\cE$}
In this paragraph, we will give lower and upper bounds for the difference
$\mathcal{E}(t) - \inf G$ as $t \to \infty$. We start with the particular case corresponding
to $G(x)=|x|^2/2$. In this case, the differential equation $(\cS)$ becomes
\begin{equation} \label{eq.linear}
\ddot{x}(t) + a(t)\, \dot{x}(t) + x(t) = 0.
\end{equation}
The next proposition precises the rate at which solutions converge to $0$.

\begin{prpstn}\label{pr.linear}
Let $a:\R_+\to \R_+$ be a non increasing map of class $\cC^2$. Assume that $\lim_{t\to\infty}
a(t)=\lim_{t\to\infty}\dot a(t)=0$ and that the map $t\mapsto \ddot a (t)+ a(t)\, \dot a (t)$
has a constant sign when $t\to \infty$. Let $x$ be a solution of the differential equation
(\ref{eq.linear}). Then there exist constants $0 < k < K < \infty$ such that for $t$ large
enough
\begin{equation}\label{alpha_est}
k e^{-\int_0^t a(s)\,ds}  \le |x(t)|^2 + |\dot{x}(t)|^2  \le K e^{-\int_0^t a(s)\,ds}\, .
\end{equation}
\end{prpstn}
\begin{proof}We eliminate the first order term in (\ref{eq.linear})
in the usual way: if we set  $A(t) =e^{\demi\,\int_0^t a(s)\,ds}$,
 then the map $y$ defined by $y(t) = A(t)x(t)$ satisfies
\begin{equation}\label{eq.equa_diff_y}
\ddot{y}(t) + \left(1- \frac{\dot{a}(t)}{2} - \frac{a(t)^2}{4}\right)y(t) = 0
\end{equation}
for every $t \geq 0$. Define the function $E:\R_+\to \R$ by
 $$E(t) = |\dot{y}(t)|^2 + \left(1-\frac{\dot{a}(t)}{2} -
\frac{a(t)^2}{4}\right) |y(t)|^2 $$ for every $t\geq 0$.
 Then $E(t)$ is non-negative for all sufficiently large $t$, and
the expression of $E(t)$ as a function of $x(t)$, $\dot x(t)$ is given by
\[E(t) = A(t)^2 \left(\left(1 - \frac{\dot{a}(t)}{2} - \frac{a(t)^2}{4}\right)|x(t)|^2 +
\left|\frac{a(t)}{2} x(t) +\dot{x}(t)\right|^2\right) \, .
\]
Therefore for sufficiently large $t$
\begin{equation}\label{eqx}
\frac{A(t)^{-2}}{2} E(t) \le |x(t)|^2 +
|\dot{x}(t)|^2 \le 2 A(t)^{-2} E(t)\, . \end{equation}
Multiplying equation (\ref{eq.equa_diff_y}) with $\dot{y}(t)$ results in
\[\dot E(t) = -\demi \left[\ddot{a}(t) +a(t)\dot{a}(t)\right]|y(t)|^2 \,.
\]
Assume now that $\ddot a (t)+ a(t)\, \dot a (t)\leq 0$ for $t$ large enough. Since $|y(t)|^2
\leq 2\, E(t)$ for sufficiently large $t$, we derive that there exists $T\geq 0$ such that
\[\forall t\geq T, \qquad 0 \le \dot E(t) \le -[\ddot{a}(t) +a(t)\dot{a}(t)]\,E(t).\]
By integrating over $[T,t]$, we obtain
$$\forall t\geq T, \qquad 0\leq \ln \frac{E(t)}{E(T)}\leq -\left[\dot a (s)+ \frac{a^2(s)}{2}
 \right]_T^t\leq \dot a (T)+ \frac{a^2(T)}{2}.$$
By setting $C=\exp\left({\dot a (T)+ \frac{a^2(T)}{2}}\right)$, we then have
\begin{equation}\label{eqy}
\forall t\geq T, \qquad E(T) \le E(t) \le C E(T).
\end{equation}
  Then estimate (\ref{alpha_est}) follows from (\ref{eqx}) and (\ref{eqy}).
If we assume that $\ddot a (t)+ a(t)\, \dot a (t)\geq 0$ for $t$ large enough, the same arguments
show that there exist $T'\geq 0$ and $C'\in (0,1)$ such that
$$\forall t\geq T', \qquad C'\,E(T) \le E(t) \le  E(T),$$
and we conclude in the same way.
\end{proof}
\begin{xmpl}\label{ex.1}
Assume that $a(t)=\frac{c}{t+1}$ for every $t\geq 0$, with $c>0$. It is immediate to check that
for $t$ large enough, $\ddot a (t)+ a(t)\, \dot a (t)\geq 0$ (resp. $\leq 0$)  if $c\leq 2$
(resp. $c\geq 2$). Therefore the assumptions of Proposition \ref{pr.linear} are satisfied and
the following estimate holds for $t$ large enough
$$ \frac{k}{t^c}\le |x(t)|^2 + |\dot{x}(t)|^2  \le \frac{K}{t^c}.$$
\end{xmpl}
\begin{xmpl}\label{ex.2}
Assume that $a(t)=\frac{1}{(t+1)^\alpha}$ for every $t\geq 0$, with $\alpha\in (0,1)$. We let the
reader check that $\ddot a (t)+ a(t)\, \dot a (t)\leq 0$ for $t$ large enough.
Therefore the assumptions of Proposition \ref{pr.linear} are satisfied and
the following estimate holds for $t$ large enough
$$ k\, e^{-t^{1-\alpha}/(1-\alpha)}\le |x(t)|^2 + |\dot{x}(t)|^2  \le
K\, e^{-t^{1-\alpha}/(1-\alpha)}.$$
\end{xmpl}
The result of Proposition \ref{pr.linear} and Examples \ref{ex.1}, \ref{ex.2}
 will serve us as a guideline in the sequel.
Let us now come back to the case of a general potential $G$. The next result provides
 a lower bound for the convergence rate of the energy
function $\cE$.  We stress the fact that there is no
convexity assumption on the function $G$ in the next statement.
\begin{prpstn}\label{prop:low_bound}
Let $a: \R_+ \to \R_+$ be a continuous map and let $G: H\to \R$ be a
function of class $\cC^1$ such that $\nabla G$ is Lipschitz continuous on the bounded
sets of $H$. If $\inf G >-\infty$, then any solution $x$ of $(\cS)$ satisfies
\begin{equation}\label{eq.mino_energie}
\forall t\geq 0, \qquad \cE(t)-\inf G \geq (\cE(0)-\inf G). \,e^{-2\, \int_0^t
a(s)\, ds}.
\end{equation}
\end{prpstn}
\begin{proof} Taking into account the expression of $\cE$ and the computation of $\dot{\cE}$,
we have
$$\dot{\cE}(t)+2\, a(t)\, (\cE(t)-\inf G)=2\, a(t)\, (G(x(t))-\inf G)\geq 0.$$
Let us multiply the above inequality by $e^{2\, \int_0^t a(s)\, ds}$, we deduce
that:
$$\forall t\geq 0, \qquad \frac{d}{dt}\left[ e^{2\, \int_0^t a(s)\, ds}\,
(\cE(t)-\min G)\right]\geq 0.$$ Formula (\ref{eq.mino_energie})
immediately follows.
\end{proof}
The next corollary gives a first result of non-convergence of the trajectories
under the condition $\int_0^{\infty} a(s)\, ds<\infty$. This hypothesis means that
the quantity $a(t)$ fastly tends to $0$ as $t\to \infty$. It is not surprising
that convergence fails under such a condition, \cf for example the extreme case $a\equiv 0$.

\begin{crllr} \label{co.non_conv:fast}
Assume that
$\int_0^{\infty} a(s)\, ds<\infty$, that the function $G$ is convex, and all the other hypotheses of Proposition \ref{prop:low_bound}. Given $(x_0, \dot x_0)\in H^2$,
consider the unique solution $x$ to the differential equation
$(\cS)$ satisfying the initial conditions $(x(0),\dot
x(0))=(x_0,\dot x_0)$.  If $(x_0,\dot x_0)\not \in \argmin G\times \{0\}$,
then the trajectory $x$ of $(\cS)$ does not converge.
\end{crllr}
\begin{proof} Let us first remark that the assumption
$(x_0,\dot x_0)\not \in \argmin G\times \{0\}$ implies that $\cE(0)>\inf G$. By taking the limit as $t\to \infty$ in inequality (\ref{eq.mino_energie}) and recalling that
$\int_0^{\infty} a(s)\, ds<\infty$, we obtain
\begin{equation}\label{eq.lim_energie}
\lim_{t\to \infty}\cE(t)-\inf G \geq (\cE(0)-\inf G). \,e^{-2\, \int_0^{\infty}
a(s)\, ds}>0.
\end{equation}
Let us now argue by contradiction and assume that there exists $\barre{x}\in H$ such that
$\lim_{t\to \infty} x(t)=\barre{x}$. From Proposition \ref{pr.stationary},
 we deduce that $\lim_{t\to \infty} \dot
x(t)=0$ and that $\nabla G(\barre{x})=0$. Since the function $G$ is convex, we infer that
 $\barre{x}\in \argmin G$. It ensues that
$\lim_{t\to \infty}\cE(t)=\min G$, which contradicts (\ref{eq.lim_energie}).
\end{proof}
The problem of convergence of the trajectories will be considered again in
section~\ref{ConvTraj}. It will be shown that condition $\int_0^{\infty} a(s)\, ds=\infty$ is
not sufficient to ensure convergence.\\

We are now going to majorize the map  $t\mapsto \cE(t)-\inf G$ as
$t\to\infty$ by some suitable quantities depending on the function $a$.
\begin{prpstn}\label{pr.speed_conv}
Let $a: \R_+ \to \R_+$ be a non increasing and differentiable map.
Let $G: H \to \R$ be a coercive function of class $\cC^1$ such that $\nabla G$ is Lipschitz continuous on the bounded sets of $H$. Assume  that $\argmin G\neq \emptyset$ and
that there exist $z\in\argmin G$ and $\theta\in \R_+$ such that condition
(\ref{eq.ineg_base}) holds.\\
%Let us denote by $x$ a solution to the differential
%equation $(\cS)$ and let $\cE$ be the associate energy function.
(i) Suppose that there exist $K_1>0$ and $t_1\geq 0$ such that $\dot{a}(t)+K_1\, a^2(t)\leq 0$,
for every $t\geq t_1$. Then, there exists $C>0$ such that, for every $t\geq t_1$,
\begin{equation}\label{eq.energy_decay1}
\cE(t)-\min G \leq C\, e^{-m\, \int_0^t a(s)\, ds},
\end{equation}
with $m= \min\left( \frac{1}{\theta+1/2},K_1\right)$.\\
(ii)  Suppose that there exist $K_2\in \left]0,\frac{1}{\theta+1/2}\right]$ and $t_2\geq 0$
such that $\dot{a}(t)+K_2\, a^2(t)\geq 0$, for every $t\geq t_2$. Then, there exists $D>0$
such that, for every $t\geq t_2$,
\begin{equation}\label{eq.energy_decay2}
\cE(t)-\min G \leq D\, a(t).
\end{equation}
\end{prpstn}
\begin{rmrk} It is immediate to see that the assumptions $\dot{a}+K_1\, a^2\leq 0$ and
$\dot{a}+K_2\, a^2\geq 0$ imply respectively that $a(t)\leq {1}/(K_1 t+c_1)$ and
$a(t)\geq {1}/(K_2 t+c_2)$,  for some $c_1$, $c_2\in \R$.
\end{rmrk}

\begin{proof} We keep the same notations as in the proof of Proposition
\ref{th.ener_L1}, in particular the expression of the map $h$ is given by
$$h(t)=\frac{a(t)}{2}\, |x(t)-z|^2+\langle \dot{x}(t), x(t)-z\rangle.$$
Let us multiply inequality (\ref{eq.relation_E_dotE}) by $e^{m\int_0^t a(s)\,
ds}$
 and integrate on the interval $[0,t]$:
\begin{equation}
e^{m\int_0^t a(s)\, ds}\,(\cE(t)-\min G)\leq \cE(0)-\min G
-\theta\,m\,\int_0^t F(s) \,\dot{h}(s)\, ds, \label{eq.exp^int}
\end{equation}
where the function $F:\R_+\to \R_+$ is defined by
$F(s)=a(s)\,e^{m\int_0^s a(u)\, du}$. The function $F$ is
differentiable and its first derivative is given by
\begin{equation}\label{eq.der_F}
\dot{F}(s)=(\dot a(s)+m\, a^2(s))\,e^{m\int_0^s a(u)\, du}.
\end{equation}
Coming back to inequality (\ref{eq.exp^int}), an integration by parts yields
$$
\int_0^t F(s) \,\dot{h}(s)\, ds=F(t)\, h(t) -F(0)\, h(0)- \int_0^t \dot F(s) \,{h}(s)\, ds.
$$
Recalling that $|h(t)|\leq M$ for every $t\geq 0$, we infer that
\begin{equation}\label{ineq.ipp}
\left |\int_0^t F(s) \,\dot{h}(s)\, ds\right|\leq M\,\left[F(t) + F(0)+ \int_0^t |\dot F(s)|\, ds
\right].
\end{equation}
We now distinguish between the cases (i) and (ii), where the assumptions allow to determine the sign of $\dot F$.
\\

\noindent (i) First assume that there exist $K_1>0$ and $t_1\geq 0$ such that
$\dot{a}(t)+K_1\, a^2(t)\leq 0$,
for every $t\geq t_1$. Let us take $m= \min\left( \frac{1}{\theta+1/2},K_1\right)$ throughout the
 proof of (i). Since $m\leq K_1$, we have $\dot a(t)+m\, a^2(t)\leq \dot a(t)+K_1\,a^2(t)\leq 0$ for every
$t\geq t_1$. It ensues from (\ref{eq.der_F})
that $\dot F(t)\leq 0$ for every $t\geq t_1$. Hence we derive from
(\ref{ineq.ipp}) that, for every $t\geq t_1$,
\begin{eqnarray*}
\left |\int_0^t F(s) \,\dot{h}(s)\, ds\right|&\leq&
M\,\left[F(t) + F(0)+ \int_0^{t_1} |\dot F(s)|\, ds+ \int_{t_1}^t |\dot F(s)|\, ds\right]\\
&=&M\,\left[F(0)+ \int_0^{t_1} |\dot F(s)|\, ds+ F(t_1)\right]=C_1.
\end{eqnarray*}
In view of (\ref{eq.exp^int}), we deduce that, for every $t\geq t_1$,
$$e^{m\int_0^t a(s)\, ds}\,(\cE(t)-\min G)\leq \cE(0)-\min G+\,\theta\, m\,C_1=C.$$
Inequality (\ref{eq.energy_decay1}) immediately follows.\\

\noindent
(ii) Now assume that there exist $K_2\in \left]0,\frac{1}{\theta+1/2}\right]$ and $t_2\geq 0$
such that $\dot{a}(t)+K_2\, a^2(t)\geq 0$, for every $t\geq t_2$. Take any $m\in
\left[K_2,\frac{1}{\theta+1/2}\right]$. Since $m\geq K_2$, we have $\dot a(t)+m\, a^2(t)\geq \dot a(t)+K_2\,a^2(t)\geq 0$ for every
$t\geq t_2$. It ensues from (\ref{eq.der_F})
that $\dot F(t)\geq 0$ for every $t\geq t_2$. Hence we derive from
(\ref{ineq.ipp}) that, for every $t\geq t_2$,
\begin{eqnarray*}
\left |\int_0^t F(s) \,\dot{h}(s)\, ds\right|&\leq&
M\,\left[F(t) + F(0)+ \int_0^{t_2} |\dot F(s)|\, ds+ \int_{t_2}^t |\dot F(s)|\, ds\right]\\
&=&M\,\left[2\, F(t)+F(0)-F(t_2)+ \int_0^{t_2} |\dot F(s)|\, ds\right]=2\,M\,F(t)+C_2.
\end{eqnarray*}
In view of (\ref{eq.exp^int}), we deduce that, for every $t\geq t_2$,
$$e^{m\int_0^t a(s)\, ds}\,(\cE(t)-\min G)\leq \cE(0)-\min G+\,\theta\, m\,C_2
+2\,\theta\, m\,M\,F(t).$$  We then
infer the existence of $C_3>0$ such that
$$e^{m\int_0^t a(s)\, ds}\,(\cE(t)-\min G)\leq C_3\, F(t),$$
which finally implies that $\cE(t)-\min G \leq C_3\, a(t)$, for every $t\geq t_2$.
\end{proof}

Let us now comment on the results given by Proposition \ref{pr.speed_conv}.
 Assume that $G$ is identically equal to $0$. In this case, condition
(\ref{eq.ineg_base}) trivially holds with $\theta=0$. If the map $a$ satisfies $\dot a +K_1\,
a^2\leq 0$ for some $K_1\geq 2$, Proposition \ref{pr.speed_conv} (i) shows that
$$\demi |\dot x(t)|^2=\cE(t)\leq C\, e^{-2\,\int_0^t a(s)\, ds}.$$
A direct computation shows that $\dot x(t)=\dot x_0\,e^{-\int_0^t
a(s)\, ds}$, so that the estimate given by Proposition
\ref{pr.speed_conv} (i) is optimal in this case. Now assume that $G$ is
given by $G(x)=|x|^2/2$ for every $x\in H$. In this case, condition
(\ref{eq.ineg_base}) holds with $\theta=1/2$ and $z=0$. Suppose that
$a(t) = \frac{c}{t+1}$ for $c \in (0,1]$. The map $a$ satisfies the
inequality $\dot a + K_1\, a^2\leq 0$ for every $K_1\in (0,\frac{1}{c}]$.
 Proposition \ref{pr.speed_conv} (i) then
 shows that $\cE(t)\leq C/(t+1)^c$ for every $t\geq 0$. Example \ref{ex.1} allows to check
that this estimate is optimal.
Suppose now that $a(t) =\frac{1}{(t+1)^\alpha}$ for some $\alpha\in (0,1)$. In this case, the map $a$
satisfies the inequality $\dot a + K_2\, a^2\geq 0$ for any $K_2>0$. Proposition \ref{pr.speed_conv} (ii) then gives the estimate  $\cE(t) \le \frac{D}{(t+1)^\alpha}$ while in fact $\cE(t) \le K\,
e^{-t^{1-\alpha}/(1-\alpha)}$ by the result of Example \ref{ex.2}.

\section{Convergence of the trajectory. Convex case}\label{ConvTraj}
Throughout this section, we are going to investigate the question of the convergence of the
trajectories associated to $(\cS)$. A first result in this direction
is given by Corollary \ref{co.non_conv:fast},
 which states that the trajectories of $(\cS)$ are not convergent
under the condition $\int_0^{\infty} a(s)\, ds<\infty$ (except for stationary solutions).
Let us now consider the particular case $G\equiv 0$. The differential
equation $(\cS)$ then becomes $\ddot x (t)+ a(t)\, \dot x (t)=0$ and
a double integration immediately shows that its solution is given by:
$$ x(t)= x(0)+\dot x(0)\,\int_{0}^{t} e^{-\int_{0}^s a(u)\, du}ds.$$
It ensues that, when $G\equiv 0$, the solution $x$ converges if and only if the
quantity $\int_0^{\infty} e^{-\int_{0}^s a(u)\, du}ds$ is finite.
Therefore it is natural to ask whether for a general potential $G$, the trajectory $x$
is convergent under the condition $\int_0^{\infty} e^{-\int_0^s a(u)\, du}ds<\infty$.
 The answer is quite complex and we will start our analysis
in the one-dimensional case.

\subsection{One-dimensional case}
First, we give a general result of non-convergence of the trajectories under the condition
\begin{equation}\label{cond_non-conv}
\int_0^{\infty} e^{-\int_0^t a(s)\, ds}dt=\infty.
\end{equation}
 Note that it is automatically satisfied if $\int_0^{\infty} a(s)\, ds<\infty$.
Condition (\ref{cond_non-conv}) expresses that the
 parametrization $a$  tends to zero rather rapidly. For example, assume that
 the map $a$ is of the form $a (t) =c/(t+1)^\gamma$,
 with $\gamma$, $c\geq 0$. It is immediate to check that
 condition (\ref{cond_non-conv}) is satisfied if and only if
$(\gamma,c)\in (1,\infty)\times \R_+$ or $(\gamma,c)\in\{1\}\times[0,1]$. Let us now state
a preliminary result.
\begin{lmm}\label{le.ineg_diff}
 Let $a:\R_+\to \R_+$ be a continuous map such that
 $\int_0^{\infty} e^{-\int_0^t a(s)\, ds}dt=~\infty$.\\
(i) Suppose that the map $p \in \cC^2([t_0,\infty),\R)$  satisfies
\begin{equation}\label{eq.lmm_ineg}
\forall t\geq t_0,\quad \ddot p (t) +a(t) \, \dot p (t) \leq 0.
\end{equation}
Then, we have either $\lim_{t\to \infty} p (t)=-\infty$ or $\,\dot p (t) \geq 0$ for every
$t\in [t_0,\infty[$.\\
(ii) Assume that the map $p$ satisfies moreover:
\begin{equation}\label{eq.lmm_eg}
\forall t\geq t_0,\quad \ddot p (t) +a(t) \, \dot p (t) = 0.
\end{equation}
Then, either $\lim_{t\to \infty}|p (t)|=\infty$ or $\,p(t)=p(t_0)$ for every
$t\in [t_0,\infty[$.
\end{lmm}

\begin{proof} (i) Assume that there exists $t_1\in [t_0,\infty[$ such that $\dot p (t_1)<0$.
Let us multiply inequality (\ref{eq.lmm_ineg}) by  $e^{\int_{t_0}^t
a(s)\, ds}$ to obtain
$$\frac{d}{dt}\left[ e^{\int_{t_0}^t a(s)\, ds} \dot p (t)\right]\leq 0.$$
Let us integrate the above inequality on the interval $[t_1,t]$, for $t\geq t_1$
$$\dot p (t)\leq \dot p (t_1)\,  e^{-\int_{t_1}^t a(s)\, ds}.$$
By integrating again, we find
$$\forall t\geq t_1, \quad p(t)\leq p(t_1)+\dot p(t_1)\,\int_{t_1}^{t}
 e^{-\int_{t_1}^s a(u)\, du}ds.$$
Recalling that $\dot p(t_1)<0$ and that $\int_{0}^{\infty}
e^{-\int_{0}^s a(u)\, du}ds
=\infty$, we conclude that $\lim_{t\to \infty} p (t)=-\infty$.\\
(ii) Assume now that the map $p$ satisfies equality (\ref{eq.lmm_eg}) and that there exists
$t_1\in [t_0,\infty[$ such that $\dot p (t_1)\neq 0$. The same computation as above shows that
$$\forall t\geq t_1, \quad p(t)= p(t_1)+\dot p(t_1)\,\int_{t_1}^{t}
 e^{-\int_{t_1}^s a(u)\, du}ds.$$
Since $\dot p(t_1)\neq 0$ and the integral $\int_{0}^{\infty}
e^{-\int_{0}^s a(u)\, du}ds$ is divergent, we conclude that
$\lim_{t\to \infty} p (t)=\pm\infty$ (depending on the sign of
$\dot p(t_1)$).
\end{proof}

Lemma \ref{le.ineg_diff} is a crucial tool in the proof of the following non-convergence result.
\begin{prpstn}\label{pr.nonconvergence}
Let $G:\R\to \R$ be a convex function of class $\cC^1$ such that $G'$ is
Lipschitz continuous on the bounded sets of $\R$.  Assume that
$\argmin G=[\alpha,\beta]$, for some $\alpha$, $\beta\in \R$ such that
$\alpha< \beta$. Let $a:\R_+\to
\R_+$ be a continuous map such that $\int_0^{\infty} e^{-\int_0^t a(s)\, ds}dt=\infty$.
 Given $(x_0, \dot x_0)\in \R^2$, consider the unique solution $x$ to
 the differential equation $(\cS)$ satisfying the initial conditions
$(x(0),\dot x(0))=(x_0,\dot x_0)$.  If $(x_0,\dot x_0)\not \in  [\alpha,\beta]\times \{0\}$,
the $\omega$- limit set $\omega(x_0,\dot x_0)$ contains $[\alpha,\beta]$, hence the trajectory
 $x$ does not converge.
\end{prpstn}

\begin{proof}  Let us assume that
\begin{equation}\label{eq.hyp}
\exists t_0\geq 0, \quad \forall t\geq t_0, \quad x(t)\geq \alpha
\end{equation}
and let us prove that this leads to a contradiction. First of all,
assertion (\ref{eq.hyp}) implies that $G'(x(t))\geq 0$ for every
$t\geq t_0$, which in view of $(\cS)$ entails that
$$\forall t\geq t_0, \quad \ddot x (t)+a(t)\, \dot x(t)\leq 0.$$
Since the map $x$ is bounded, we deduce from Lemma
\ref{le.ineg_diff}(i) that $\dot x(t)\geq 0$ for every $t\geq t_0$, hence $\barre{x}:=\lim_{t\to
\infty} x(t)$ exists. From Proposition \ref{pr.stationary}, we
have $G'(\barre{x})=0$, \ie $\barre{x} \in [\alpha,\beta]$. It ensues that
 $x(t)\in [\alpha,\beta]$ for
every $t\geq t_0$. Hence we have $G'(x(t))=0$ for every $t\geq t_0$
and we infer that
$$\forall t\geq t_0, \quad \ddot x (t)+a(t)\, \dot x(t)= 0.$$
Since  $x([t_0,\infty[)\subset [\alpha,\beta]$, we derive from Lemma \ref{le.ineg_diff}(ii) that $x(t)=x(t_0)$ for
every $t\geq t_0$. Thus, it follows by backward uniqueness that we have a stationary solution,
which contradicts the assumption $(x_0,\dot x_0)\not \in [\alpha,\beta]\times \{0\}$. Hence, we
deduce that assertion (\ref{eq.hyp}) is false, so that we can build a sequence $(t_n)$ tending
to $\infty$ such that $x(t_n)<\alpha$.  In a symmetric way, we can construct a sequence $(u_n)$
tending to $\infty$ such that $x(u_n)>\beta$. Recalling that the $\omega$-limit set $\omega(x_0,\dot x_0)$
is connected, we conclude that $\omega(x_0,\dot x_0)\supset [\alpha,\beta]$.
\end{proof}

We can now wonder if the converse assertion is true: do the trajectories
$x$ of $(\cS)$ converge under the condition $\int_0^{\infty} e^{-\int_0^s a(u)\, du}ds<\infty$?
When the coefficient $a(t)$ is constant and positive,
 the trajectories of the so-called ``Heavy Ball with Friction'' 
system are known to be convergent, see for example \cite{AttGouRed}. The question
 is more delicate in the case of an asymptotically vanishing map $a$.
 We mention below a first positive result when the map $a$ is of the form
 $a(t)=\frac{c}{(t+1)^\gamma}$ with $c$, $\gamma>0$.
\begin{prpstn}
Let $G:\R\to \R$ be a convex function of class $\cC^1$ such that $G'$ is
Lipschitz continuous on the bounded sets of $\R$.  Assume that
$\argmin G=[\alpha,\beta]$ with $\alpha \leq \beta$ and that there exists $\delta>0$ such that
$$\forall \xi\in (-\infty,\alpha], \quad G'(\xi)\leq  \delta\, (\xi-\alpha)\quad \mbox{ and }
\quad \forall \xi\in [\beta, \infty), \quad G'(\xi)\geq  \delta\, (\xi-\beta).$$
Given $c$, $\gamma>0$, let $a:\R_+\to \R_+$ be the map defined by $a(t)=\frac{c}{(t+1)^\gamma}$
for every $t\geq 0$. If $\gamma \in (0,1)$ or if $\gamma=1$ and $c>1$,  then for any solution
$x$ of $(\cS)$, $\lim_{t\to\infty} x(t)$ exists.
\end{prpstn}
We omit the proof of this result since it is rather technical and it will be developped more widely
in a future paper.

\subsection{Multi-dimensional case}
Our purpose now is to extend the result of Proposition \ref{pr.nonconvergence}
to the case of a dimension greater than one.
 The situation is much more complicated since we have to take into account the geometry of the set
$\argmin G$. In the sequel, we will assume that the gradient of $G$ satisfies the following
 condition:
\begin{equation}\label{lim_dir_grad}
\forall \barre{x}\in \bd(\argmin G), \qquad\lim_{x\to \barre{x},\,\, x\not \in \argmin G}
\frac{\nabla G(x)}{|\nabla G(x)|}\quad \mbox{ exists.}
\end{equation}
If $G$ is a convex function of a single variable, \ie $H=\R$, this condition is satisfied when
 the set $\argmin G$ is not reduced to a singleton.
Before stating the main result of non-convergence for the trajectories of $(\cS)$,
 let us first recall some basic notions of convex analysis.
The polar cone $K^\ast$ of a cone $K\subset H$ is defined by
$$K^\ast=\{y\in H|\quad \forall x\in K, \quad \langle x,y\rangle\leq 0\}.$$
Let $S\subset H$ be a closed convex set and let $\barre{x}\in S$. The normal cone $N_S(\barre{x})$
and the tangent cone $T_S(\barre{x})$ are respectively defined by
\begin{eqnarray*}
N_S(\barre{x})&=&\{\xi\in H|\quad \forall x\in S, \quad \langle \xi, x-\barre{x}\rangle\leq 0\}\\
T_S(\barre{x})&=&{\rm cl}\left[\cup_{\lambda>0}\lambda\, (S-\barre{x})\right].
\end{eqnarray*}
The convex cones $N_S(\barre{x})$ and $T_S(\barre{x})$ are polar to each other,
 \ie $N_S(\barre{x})=[T_S(\barre{x})]^\ast$ and
$T_S(\barre{x})=[N_S(\barre{x})]^\ast$. Fur further details relative to convex
analysis, the reader is referred to Rockafellar's book \cite{rockafellar}.\\

The next lemma brings to light a geometrical property that will be crucial in the proof of the next theorem.

\begin{lmm}\label{le.lim_dir_grad}
Let $G:H\to \R$ be a convex function of class $\cC^1$ such that $S=\argmin G\neq~\emptyset$.
 Given
$\barre{x}\in \bd(S)$, assume that $d:=\lim_{x\to \barre{x},\,\, x\not \in S} \nabla G(x)/
|\nabla G(x)|$ exists. Then we have:\\
(i) $N_S(\barre{x})=\R_+\, d$.\footnote{When the normal cone
$N_S(\barre{x})$ is reduced to a half-line, the set $S$ is said to be
{\it smooth} at $\barre{x}$. Hence item (i) shows that the existence
of $\lim_{x\to \barre{x},\,\, x\not \in S} \nabla G(x)/|\nabla G(x)|$ implies
smoothness of the set
$S$ at $\barre{x}$.}\\
(ii) There exists a neighborhood $V$ of $\barre{x}$, a closed convex
cone $K\subset H$ along with a positive real $\eta>0$ such that:
$$\forall x\in V, \quad \nabla G(x)\in K\quad \mbox{ and } \quad
-K\cap \eta\, \BB\subset S-\barre{x}.\leqno (\cC_{\barre{x}})$$
\end{lmm}

\begin{proof} (i) Since the function $G$ is convex and since $\nabla G(\barre{x})=0$, we
have
\begin{equation}\label{eq.monoton}
\forall x\in H, \quad \langle \nabla G(x), x-\barre{x}\rangle \geq 0.
\end{equation}
Let $v\not \in T_S(\barre{x})$ and take the vector $x$ equal to
$\barre{x}+tv$, for some $t>0$. Remark that, since $v\not \in
T_S(\barre{x})$, we have $\barre{x}+tv\not \in S$ and hence
$\nabla G(\barre{x}+tv)\neq 0$ for every $t>0$. From (\ref{eq.monoton}), we
derive that
$$\forall t>0, \quad \left\langle \frac{\nabla G(\barre{x}+tv)}{|\nabla G(\barre{x}+tv)|},
 v\right\rangle \geq 0.$$
Taking the limit as $t\to 0^+$, we infer that $\langle d,v\rangle
\geq 0$. Since this is true for every $v\not \in T_S(\barre{x})$, we
conclude that $-d\in [H\setminus T_S(\barre{x})]^\ast$. Let us denote
by $\cH_{x,\leq}$ (resp. $\cH_{x,>}$) the closed (resp. open)
hyperplane defined by $\cH_{x,\leq}=\{y\in H,\quad \langle
x,y\rangle \leq 0\}$ (resp. $\cH_{x,>}=\{y\in H,\quad \langle
x,y\rangle > 0\}$). The polarity relation
$T_S(\barre{x})=[N_S(\barre{x})]^\ast$ can be equivalently rewritten as
$\displaystyle T_S(\barre{x})=\bigcap_{x\in N_S(\barre{x})}\cH_{x,\leq}$.
Then it follows that
$$-d\in [H\setminus T_S(\barre{x})]^\ast=\Big[\bigcup_{x\in N_S(\barre{x})}\cH_{x,>}\Big]^\ast
=\bigcap_{x\in N_S(\barre{x})}\cH_{x,>}^\ast=\bigcap_{x\in
N_S(\barre{x})}\R_-\,x.$$
If the cone $N_S(\barre{x})$ is not reduced to a half-line, the above intersection equals $\{0\}$,
which contradicts the fact that $|d|=1$. Hence the cone $N_S(\barre{x})$ is equal to a half-line
and the above inclusion shows that $N_S(\barre{x})=\R_+.d$.\\
(ii) Since $\lim_{x\to \barre{x},\,\, x\not \in S}
\nabla G(x)/|\nabla G(x)|=d$, there exists a neighbourhood $V$ of $\barre{x}$ such
that
$$\forall x\in V, \quad x\not \in S \quad \Longrightarrow\quad \left\langle d, \frac{\nabla G(x)}
{|\nabla G(x)|} \right\rangle\geq \demi.$$ Let us define the cone
$K=\{0\}\cup\{v\in H,\quad \langle d, \frac{v}{|v|} \rangle\geq \demi\}$.
 It is clear that $K$ is a closed convex cone and that
$\nabla G(x)\in K$ for every $x\in V$. On the other hand, since
$N_S(\barre{x})=\R_+.d$, we have
$$\liminf_{v\to 0, \, \, v\not \in S-\barre{x}}\langle d, \frac{v}{|v|} \rangle\geq 0.$$
Therefore, there exists $\eta>0$ such that
$$v\in \eta\, \BB\quad \mbox{ and }\quad v\not \in S-\barre{x}  \quad \Longrightarrow\quad
\langle d, \frac{v}{|v|} \rangle >-\demi.$$ Since for every
$v\in -K\setminus \{0\}$, $\langle d, \frac{v}{|v|} \rangle \leq -\demi$,
 we deduce that $-K\cap \eta\, \BB\subset S-\barre{x}$,
which achieves the proof of property $(\cC_{\barre{x}})$.
\end{proof}

Let us now state the general result of non-convergence for the trajectories of $(\cS)$
under the condition $\int_0^{\infty} e^{-\int_0^t a(s)\, ds}dt=\infty.$
\begin{thrm}\label{th.non-conv}
Let $G:H\to \R$ be a convex coercive function of class $\cC^1$ such that $\nabla G$ is Lipschitz
 continuous on the bounded sets of $H$. Assume moreover that the geometric property (\ref{lim_dir_grad}) holds.
 Let $a:\R_+\to \R_+$ be a continuous map such that
$\int_0^{\infty} e^{-\int_0^t a(s)\, ds}dt=\infty$. Given $(x_0, \dot x_0)\in H^2$,
consider the unique solution $x$ of $(\cS)$ satisfying  $(x(0),\dot
x(0))=(x_0,\dot x_0)$.  If $(x_0,\dot x_0)\not \in \argmin G\times \{0\}$,
then the trajectory $x$ of $(\cS)$ does not converge.
\end{thrm}

\begin{proof} For simplicity of notation, let us set $S=\argmin G$.
Let us prove the contraposition of the previous statement and
assume that there exists $\barre{x}\in H$ such
that $\lim_{t\to \infty} x(t)=\barre{x}$.  We must prove that in this case $x_0 \in S, \, \dot{x}_0 = 0$. From Proposition \ref{pr.stationary}, we
have $\nabla G(\barre{x})=0$, hence $\barre{x} \in S$. If $\barre{x}\in \bd(S)$,
condition  $(\cC_{\barre{x}})$ is satisfied in view of assumption (\ref{lim_dir_grad}) and
Lemma \ref{le.lim_dir_grad} (ii). On the other hand, if $\barre{x}\in \inte(S)$ condition  $(\cC_{\barre{x}})$ is trivially satisfied with $K=\{0\}$.
In both cases, we derive the existence of a closed convex cone
$K\subset H$ along with $\eta>0$ and $t_0\geq 0$ such that
\begin{equation}\label{eq.C_x}
\forall t\geq t_0,\quad \nabla G(x(t))\in K\quad \mbox{ and }\quad -K\cap
\eta\, \BB\subset S-\barre{x}.
\end{equation}
Let $v\in -K^\ast$ and take the scalar product of $(\cS)$ by the
vector $v$. Since $\langle \nabla G(x(t)),v \rangle \geq 0$ for every
$t\geq t_0$, we deduce that
$$\forall t\geq t_0,\quad \langle \ddot x(t),v\rangle+a(t)\, \langle \dot x(t),v\rangle \leq 0.
$$
Let us apply Lemma \ref{le.ineg_diff}(i) to the map $p$ defined by
$p(t)=\langle  x(t),v\rangle$. Since the trajectory $x$ is bounded,
we infer that $\langle \dot x(t),v\rangle \geq 0$ for every $t\geq
t_0$. By integrating on the interval $[t,\infty[$, we find $\langle
x(t)-\barre{x},v\rangle \leq 0$ for every $t\geq t_0$. Since this is
true for every $v\in -K^\ast$, we derive that $x(t)-\barre{x}\in
-K^{\ast\ast}$ for every $t\geq t_0$. Recalling that
$K^{\ast\ast}=K$ for every closed convex cone $K$, we conclude that
$x(t)-\barre{x}\in -K$ for every $t\geq t_0$. On the other hand, since
$\lim_{t\to \infty} x(t)=\barre{x}$, there exists $t_1\geq t_0$ such
that $x(t)-\barre{x}\in \eta\, \BB$ for every $t\geq t_1$. In view of
(\ref{eq.C_x}), we infer that $x(t)\in S$ for every $t\geq t_1$, so
that the differential equation $(\cS)$ becomes
$$\forall t\geq t_1,\quad \ddot x(t)+a(t)\, \dot x(t)=0.$$
By arguing as in the proof of  Lemma \ref{le.ineg_diff}(ii), we deduce that either
$\lim_{t\to \infty} |x(t)|=\infty$ or $x(t)=x(t_1)$ for every $t\geq t_1$. Since the map $x$
converges toward $\barre{x}$, the first eventuality does not hold. It follows by backward
 uniqueness that we have a stationary solution, $x(t) = x_0$ for all $t$, which must therefore satisfy
$(x_0,\dot x_0) \in S\times \{0\}$.
\end{proof}

\section{The case of a Non-Convex Potential}\label{sec:NonConv}
In this section, we discuss the case where $G$ is defined on
$\R^n$ and has multiple critical points, but does not necessarily satisfy
condition (\ref{eq.ineg_base}). Instead, we will assume that

(a) $G$ has finitely many critical points $x_1, \, x_2, \, \dots, x_N$.

(b) $G$ attains different values on them, i.e. we can order them such that
\[\lambda_1 = G(x_1) < \lambda_2 = G(x_2) < \lambda_3 = G(x_3) < \dots \, < \lambda_N = G(x_N). \]
This is the "generic case". We will also use the assumption

(c) $a:\R_+ \to \R_+$ is non increasing such that $\int_0^\infty a(s) ds = \infty$ .

\smallskip
Our first result shows that in this case, for each solution there exists exactly one critical point that is visited for arbitrarily long times.

\medskip
\begin{prpstn} \label{pr.accumulation}
Let $G:\R^n\to \R$ be a coercive function of class $\cC^1$ such that
$\nabla G$ is Lipschitz continuous
on the bounded sets of $\R^n$. If the assumptions (a)-(c) above are satisfied, then
 there exists a unique $x^* \in \{x_1, x_2, \dots, x_N\}$ such that
\[\lim_{t \to \infty}\cE(t) = G(x^*) \, .\]
Also, for all $T\geq 0$
\[\liminf_{t \to \infty} \sup_{s \in [t,t+T]} |x(s)-x^*| = 0. \]
Moreover  $x^*$ is the only point of accumulation that is visited for arbitrarily long time
intervals. If $x^*$ is a local minimum of $G$, then in fact
\[\lim_{t \to \infty} |x(t)-x^*| = 0. \]
\end{prpstn}

\begin{proof} From Proposition \ref{pr.slowdown}, there are arbitrary long time intervals
$[S,T]$ where $|\dot{x}(t)|$ is arbitrarily small.
Since $t\mapsto \cE(t)=\demi |\dot x (t)|^2+G(x(t))$ is decreasing  and bounded
below, its limit exists. If this is not equal to $G(x_j)$ for some $j\in\{1,\hdots,N\}$, then
there exists an interval  $[a,b]$  that does not contain any critical value of $G$, such that
$G(x(t)) \in [a,b]$ for every $t\in [S,T]$. Now, we can find $\delta>0$ such that $|g(\xi)| \ge \delta$ on $G^{-1}([a,b])$, and thus $|g(x(t))| \ge \delta$ for every $t\in [S,T]$. However,
this contradicts Proposition \ref{pr.g.estimate}. Thus, there exists
$x^* \in \{x_1, x_2, \dots, x_N\}$ such that $\lim_{t \to \infty}\cE(t) = G(x^*)$.
Therefore, $|x(t) - x^*|$ becomes arbitrarily small on arbitrarily long
intervals by Proposition \ref{pr.g.estimate}. To show that there is no other
point $x^{**}$ such that
\[\liminf_{t \to \infty} \sup_{s \in [t,t+T]} |x(s)-x^{**}| = 0
\, ,\]
\medskip
firstly suppose that such a point exists with $G(x^{**}) = G(x^*)$. Note that then
$g(x^{**}) \ne 0$. On the other hand, $|\dot x(t)|$ must become arbitrarily small on these same
intervals, since $\lim_{t\to \infty}\cE(t) = G(x^*) = G(x^{**})$, and therefore also $g(x(t))$
can be made arbitrarily small on such intervals. But this impossible, since
$g(x^{**}) \ne 0$.

Next, suppose that $G(x^{**}) < G(x^*)$. In this case, we can find arbitrarily
long intervals $[S,T]$ where $|x(t) - x^{**}|$ is small and therefore $G(x(t))
< G(x^*) - \delta$ for some $\delta>0$. Then $ |\dot{x}(t)|^2 \ge 2 \delta$ on
these intervals. However, the map $\ddot x$ is uniformly bounded on any such interval.
By applying Landau's inequality on the interval $[S,T]$,
\[\|\dot{x}\|_{\infty} \le 2 \sqrt{\|x - x^{**}\|_{\infty} \|\ddot{x}\|_{\infty}}\, ,
\]
we obtain a contradiction. Hence no such point $x^{**}$ exists.

\medskip
Since $\lim_{t\to \infty}\cE(t) = G(x^*)$, the solution  $x(t)$ must eventually
enter and remain in the connected component of $G^{-1}((-\infty,G(x^*)+\eps))$ that contains $x^*$, for any $\eps > 0$.
If $^*$ is a local minimum of $G$ (which is strict by assumption), then the intersection of these neighborhoods of $x^*$ is just
$\{x^*\}$; hence $\lim_{t \to \infty} x(t) = x^*$.
\end{proof}

We now give a result that shows that the density of the times $t\in \R_+$ when $x(t)$
is near the critical point $x^*$ approaches 1. This is comparable to a result on "convergence in probability".
\medskip
\begin{thrm} \label{th.density}
In addition to the hypotheses of Proposition \ref{pr.accumulation}, we assume that $\dot a \in L^1_{loc}(\R_+)$ and that there exists $c>0$ such that $a(t) \ge \frac{c}{t+1}$ for every $t\geq 0$. Then there is a unique stationary point $x^* \in \{x_1,\, \dots, x_N\}$ of $G$ such that for any $\eps >0$
\begin{equation} \lim_{T \to \infty} T^{-1} \big| \{ t \le T | |x(t) - x^*| > \eps \} \big| = 0 \, ,
\end{equation}
where $\big| A \big|$ denotes the one-dimensional Lebesgue measure of a measurable set $A \subset \R$.
\end{thrm}

\begin{proof} We first recall estimate (\ref{energy3}) which implies
\begin{equation}\label{dotx.est}
\int_0^\infty \frac{1}{t+1} |\dot{x}(t)|^2 dt < \infty \, .
\end{equation}
We are now going to show that
\begin{equation}\label{g.est}
\int_0^\infty \frac{1}{t+1} |g(x(t))|^2 dt  < \infty.
\end{equation}
 We write $b(t) = \frac{1}{t+1}$. Form the scalar product of $(\cS)$ with $b(t) g(x(t))$ and integrate over $[0,T]$, for some $T>0$. The result is the identity
$$\int_0^T b(t)|g(x(t))|^2 dt = - \int_0^T a(t)b(t) \langle \dot{x}(t), g(x(t)) \rangle dt
 -  \int_0^T b(t) \langle \ddot{x}(t), g(x(t)) \rangle dt  \, .$$
Integrating by parts, the first integral on the right hand side becomes
\[- \int_0^T a(t)b(t) \frac{d}{dt} G(x(t)) dt = \left( a(t)b(t) G(x(t)) \right) \big|_{t=T}^{t=0}  + \int_0^T \frac{d}{dt}\left(a(t) b(t)  \right) G(x(t)) dt \, .
\]
Since $G(x(\cdot))$ is bounded on $[0,\infty)$ and $\dot{a}b + a \dot{b}$ is integrable, this term therefore is bounded. The second integral on the right hand side of the above identity becomes after two integrations by parts
\begin{eqnarray*}
\dots &=&  \left( -b(t) \langle \dot{x}(t), g(x(t)) \rangle \right) \big|_{t=0}^{t=T} + \int_0^T \dot{b}(t) \langle \dot{x}(t), g(x(t)) \rangle dt \\
&& + \int_0^T b(t) \langle \dot{x}(t), \frac{d}{dt} g(x(t)) \rangle dt  \\
&=& \left(-b(t) \langle \dot{x}(t), g(x(t)) \rangle + \dot{b}(t) G(x(t)) \right) \big|_{t=0}^{t=T} - \int_0^T \ddot{b}(t) G(x(t))  dt \\
&& + \int_0^T b(t) \langle \dot{x}(t), Dg(x(t)) \dot{x}(t)\rangle dt
\end{eqnarray*}
where $Dg(x(t))$ is the derivative of $g$ at $x(t)$. The first two terms are both bounded due to previous estimates, and $\|Dg(x(t))\| \le M$ for all $t$ since the trajectory $x(\cdot)$ is bounded and $g = \nabla G$ is Lipschitz, uniformly on bounded sets. Therefore the last integral is bounded in magnitude by
\[\left| \int_0^T b(t) \langle \dot{x}(t), Dg(x(t)) \dot{x}(t)\rangle dt \right|
\le M \int_0^T b(t) |\dot{x}(t)|^2 dt \le \frac{M}{c} \int_0^T a(t) |\dot{x}(t)|^2 dt
\]
which remains uniformly bounded for all $T>0$.  This proves (\ref{g.est}).

\smallskip
By Proposition \ref{pr.accumulation}, $\lim_{t \to \infty} \cE(t) = G(x^*) = \lambda_i$ for some stationary point $x^* = x_i$ of $G$. Pick $\delta > 0$ such that $\lambda_{i-1} + \delta < \lambda_i< \lambda_{i+1} - \delta$, and let $T_0>0$ be so large that $\cE(t) < \lambda_i+\delta$ if $t \ge T_0$. For the remainder of the proof, we assume without loss of generality that $T_0 = 0$.
Let $\eps > 0$, then there exists $\gamma >0$ such that
\begin{equation}\label{eq.implic}
|x(t) - x^*| \ge \eps \quad \Rightarrow \quad |\dot{x}(t)| \ge \sqrt{2 \delta} \quad \text{or} \quad |g(x(t))| \ge \gamma  \,.
\end{equation}
Indeed, assume that  $|x(t) - x^*| \ge \eps$ and that $|\dot{x}(t)| < \sqrt{2 \delta}$.
First we have
\[G(x(t)) = \cE(t) - \frac{1}{2}|\dot{x}(t)|^2 \ge  \lambda_i - \delta \,.
\]
Now set
\[\gamma = \min \{|g(\xi)| \, \big| \, |\xi - x^*| \ge \epsilon, \, \lambda_i - \delta \le G(\xi) \le \lambda_i + \delta \, \} \,.\]
The quantity $\gamma$ is positive since there are no critical points of $G$ in the compact region over which $|g(\cdot)|$ is minimized. Hence assertion (\ref{eq.implic}) is proved. Therefore,
we deduce that
\begin{eqnarray*} \big| \{t \big||x(t) - x^*| \ge \epsilon\}\big| &\le& \big|\{t \big| |\dot{x}(t)| \ge \sqrt{2 \delta} \quad \text{or} \quad |g(x(t))| \ge \gamma \}\big|\\
&\le& \big|\{t \big||\dot{x}(t)| \ge \sqrt{2 \delta}\}\big| + \big|\{t \big| |g(x(t))| \ge \gamma \}\big| \, .
\end{eqnarray*}
By estimates (\ref{dotx.est}) and (\ref{g.est}) combined with Lemma \ref{lmm.density},
we derive that
$$\hspace{-7cm}
\lim_{T \to \infty} T^{-1} \big| \{ t \le T \big| |x(t) - x^*| \ge \eps \} \big|\le$$
$$\hspace{1cm}\lim_{T \to \infty} T^{-1} \big| \{ t \le T \big| |\dot{x}(t)| \ge \sqrt{2 \delta}\} \big|+\lim_{T \to \infty} T^{-1} \big| \{ t \le T \big| |g(x(t))| \ge \gamma \} \big|=0. $$
The theorem has been proved.
\end{proof}

As a consequence, a Cesaro average of the solution $x(\cdot)$ converges to the
 critical point $x^*$.

\begin{crllr}\label{crl.ergodic} Under the conditions of Theorem \ref{th.density}, all solutions of $(\cS)$ satisfy
\[ \lim_{T \to \infty} \frac{1}{T} \int_0^T x(t) dt = x^* \]
for some critical point $x^*$ of $G$.
\end{crllr}

\begin{proof} Let $x^*$ be the stationary point identified in Theorem \ref{th.density}. Given $\eps > 0$ and $T >0$, we have
\[ \left| \frac{1}{T}\int_0^T x(t) dt -  x^* \right| \le \frac{1}{T} \int_0^T |x(t) - x^*| dt \le \eps +
\frac{M}{T} \big|\{t\le T | |x(t) - x^*| \ge \eps \} \big|\]
where $M = \sup_{t \ge 0} |x(t) - x^*| < \infty$. The $\limsup$ of the right hand side is no larger than $\eps$, which proves the corollary.
\end{proof}

In the same direction, a result of convergence of the ergodic mean has been obtained by Brezis \cite{brezis} for trajectories associated to semigroups of nonlinear contractions.\\

Let us now establish a result that was useful in the proof of Theorem \ref{th.density}.

\medskip
\begin{lmm} \label{lmm.density} Let $w : [0,\infty) \to \R$ be measurable, bounded, and non-negative. If
\[ \int_0^\infty \frac{w(t)}{t+1} dt < \infty \, ,\]
then for all $\eps > 0$
\[\lim_{T \to \infty} T^{-1} \big| \{ t \le T | w(t) > \eps \} \big| = 0 \, .
\]
\end{lmm}

\begin{proof}Let $\eps > 0$. For $0 \le S < T < \infty$, set $A(\eps,S,T) = \{ S \le t \le T | w(t) > \eps \}$. It is clearly sufficient to show that for all $\delta \in (0,1)$ we can find $S > 0$ such that
\[\limsup_{T \to \infty} (T+1)^{-1} \big| A(\eps,S,T) \big| \le \delta \, .
\]
For this purpose, define $\rho = - \log(1 - \delta)$, and choose $S > 0$ large enough such that $\int_S^\infty \frac{w(t)}{t+1}dt \le \eps \rho$. Then for $T>S$,
\[
\eps \rho \ge \int_S^T \frac{w(t)}{t+1}dt \ge \int_{A(\eps,S,T)} \frac{\eps}{t+1} dt. \]
Clearly the integral on the right is minimized if $A(\eps,S,T) = [T-  \big| A(\eps,S,T) \big|,T]$, and therefore
\[\rho \ge \int_{T-  \big| A(\eps,S,T) \big|}^T \frac{1}{t+1} dt =
 - \log \left(1 - \frac{\big| A(\eps,S,T) \big|}{T+1} \right)\, .\]
By the choice of $\rho$, this inequality becomes $\frac{\left| A(\eps,S,T) \right|}{T+1} \le \delta \, .$
\end{proof}
We stress the fact that the above proof mainly relies on a variant of Markov's inequality.\\

The main remaining open question is whether  $\lim_{t \to \infty}x(t)$ exists
if $a(t) \to 0$ as $t \to \infty$. There is a unique stationary point of $G$ that is visited for long times. If this point is a local minimum of $G$, the trajectory will converge to it by Proposition \ref{pr.accumulation}. If the point is a local maximum, it appears possible to adapt the arguments from the next section to show that again convergence holds. The difficulty is that in more than one dimension, the stationary point that is visited for long times
may be a saddle point of $G$. Then it is possible that the solution visits other regions
of $\R^n$ intermittently for finite amounts of time, infinitely often, spending longer
 and longer periods of time near the saddle point in between. In one dimension,
such a behavior cannot occur, since solutions either get trapped near local minima of $G$,
or if they visit local maxima of $G$ without converging to them, they must leave their neighborhood rapidly. \\

To end this section, we show that if $a(\cdot)$ is bounded away from 0 (e.g. if $a(\cdot)$ is a positive constant), solutions of ($\cS$) always converge.
\begin{prpstn}\label{pr.bounded_a}
Let $a:\R_+ \to \R_+$ be a non increasing map such that $a(t) \ge a_0$ for all
$t\geq 0$, with some $a_0>0$. Let
$G:\R^n\to \R$ be a coercive function of class $\cC^1$, such that
$\nabla G$ is Lipschitz continuous on the bounded sets of $\R^n$. Assume that $G$ has finitely
many critical points $x_1, \, x_2, \, \dots, x_N$. Then, for any solution $x$
to the differential equation $(\cS)$, there exists $x^{*} \in \{x_1, \, x_2, \, \dots, \, x_N \}$
such that $\lim_{t \to \infty} x(t) = x^{*}.$
\end{prpstn}
\begin{proof} The assumption $a(t)\ge a_0>0$ implies that $\int_0^\infty |\dot{x}(t)|^2 dt < \infty$
and hence $\lim_{t \to \infty}\dot{x}(t)=0$, since the map $\ddot{x}$ is uniformly bounded.
 From Proposition \ref{pr.g.estimate}, we derive that $\lim_{t \to \infty} g(x(t))=0$.
Since the set of zeroes of $g$ is discrete, this implies that $x(t)$ converges to one
of the critical points of $G$.
\end{proof}

\section{The One-Dimensional Case} \label{sec:OneDim}
Let us consider the equation $(\cS)$ in the one-dimensional case.
The derivative $\dot x$ changes sign either finitely many times or infinitely many times. In
the first case, solutions must have a limit, while the second case can occur
either if the solution approaches a limit or if the $\omega$-limit of the
trajectory is a non-empty interval. We shall give conditions that exclude this
last possibility. Rather, trajectories always have a limit, and moreover
solutions oscillate infinitely if and only if this limit is a local minimum of
$G$. We show further that in this case the set of initial conditions for which
solutions converge to a local minimum is open and dense.

\medskip
To describe the behavior of the trajectories more precisely, let us write
$$w(t) = \cE(t) = G(x(t)) + \frac{1}{2} |\dot{x}(t)|^2$$ and observe that
\[\dot w(t) = -a(t)|\dot{x}(t)|^2 = 2a(t)\left(G(x(t)) - w(t)\right) \]
and
\[\dot{x}(t) = \pm \sqrt{2w(t) - 2G(x(t))}\,.\]
Assume that $a(t)>0$ for every $t\geq 0$ and that the solution $x$ is not stationary.
It is obvious that $\dot w(t) <0$ for $t\geq 0$, except at times $t$ where $w(t) = G(x(t))$,
and these $t$ are precisely those times
where $\dot{x}$ changes sign.  The set $\cT=\{t\geq 0\,|\,\dot x(t)=0\}$ must be discrete.
 Indeed, if $t^*$ is an accumulation point of $\cT$, then there exists a sequence
$(t_i)$ tending toward $t^*$ such that $\dot x (t_i)=0$, hence $\dot x(t^*)=0$. By Rolle's
Theorem, there exists also a sequence $(u_i)$ tending toward $t^*$ such that $\ddot x (u_i)=0$, hence $\ddot x(t^*)=0$. Hence we would have the equality $\dot x(t^*)=\ddot x(t^*)=
G'(x(t^*))=0$ and thus $x$ would have to be a constant solution, a contradiction. Therefore
 there exists an increasing sequence $(t_n)$ tending toward $\infty$ such that
$\cT=\{t_n,\, n\in \N\}$. As $\cT$ is discrete, $w$
is strictly decreasing, hence $\lim_{t\to \infty} w(t)$ exists if the function
$G$ is bounded from below. If $\cT$ is finite, \textit{i.e.} if $\dot{x}$ changes its sign finitely
often, then $\lim_{t \to \infty} x(t) = x^*$ exists since $x$ is eventually monotone and is bounded (provided some coercivity assumption on $G$). In this case, $G'(x^*) = 0$ by Proposition
 \ref{pr.stationary}.
 However, without additional assumptions on the maps $a$ and $G$, the trajectory $x(\cdot)$
needs not converge, as Proposition \ref{pr.nonconvergence} shows.

\medskip
Before giving the main assumptions of this section, let us recall the definitions of strong
convexity and strong concavity.
\begin{definition}
Let $G:\R\to \R$ a function of class $\cC^1$ and let $x^*\in \R$. The function $G$ is said
to be strongly convex in the neighbourhood of $x^*$ if there exist $\eps$, $\delta>0$ such
that
\begin{equation}\label{eq.strong_conv}
\forall x, y\in ]x^*-\eps, x^*+\eps[, \qquad G(y)\geq G(x)+ (y-x)\, G'(x)+\, \delta\, (y-x)^2.
\end{equation}
\end{definition}
It is easy to check that the above property amounts to saying that the map $x\mapsto G(x)-\delta\, x^2$ is convex on $]x^*-\eps, x^*+\eps[$. This is also equivalent to the fact that the
map $x\mapsto G'(x)-2\,\delta\, x$ is non decreasing on $]x^*-\eps, x^*+\eps[$. When the function
$G$ is of class $\cC^2$, assertion (\ref{eq.strong_conv}) is equivalent to the inequality $G''\geq 2\, \delta$ on $]x^*-\eps, x^*+\eps[$.\\
Let us now introduce the notion of strong concavity.
\begin{definition}
Let $G:\R\to \R$ a function of class $\cC^1$ and let $x^*\in \R$. The function $G$ is said
to be strongly concave in the neighbourhood of $x^*$ if $-G$ is strongly convex
in the neighbourhood of $x^*$.
\end{definition}

We are now able to set up the framework that will be used throughout this section. The function
 $G:\R\to \R$ of class $\cC^1$ will satisfy the following assumptions, considered as the "generic" case.

\begin{itemize}
\item[] (a) $G$ has finitely many critical points $x_1 < x_k < \dots < x_N$.
\item[] (b) If $k \neq r$, then $G(x_k) \neq G(x_r)$.
\item[] (c) For all $k\in \{1,\hdots,N\}$, the function $G$ is either strongly convex or
strongly concave in the neighbourhood of $x_k$.
\end{itemize}
\medskip
Property (c) implies that the critical points of $G$ correspond either to local minima or local
maxima of $G$. Moreover, property (c) shows that near local minima $x_j$, $G$ satisfies the
inequality $G(x) \ge G(x_j) + \delta |x-x_j|^2$ and near local maxima $x_k$, we
have similarly $G(x) \le G(x_k) - \delta |x-x_k|^2$. We can now describe the
asymptotic behavior of solutions of $(\cS)$.

\medskip
\begin{thrm}\label{thm_a}
Let $a:\R_+ \to \R_+$ be a differentiable non increasing map such that
$ \lim_{t\to\infty} a(t)=0$. Assume that there exists $c>0$ such that $ a(t)\geq \frac{c}{t+1}$
for every $t\geq 0$. Let $G:\R\to \R$ be a coercive function of class $\cC^1$ such that
$G'$ is Lipschitz continuous on the bounded sets of $\R$. If $G$ satisfies the additional
assumptions (a)-(b)-(c) above, then for any solution $x$ of $(\cS)$, $\lim_{t \to \infty} x(t)$
exists. Moreover, denoting by $\cT$ the set of sign changes of $\dot x$,
the limit is a local maximum of $G$ if and only if the set $\cT$ is finite,
and it is a local minimum of $G$ if and only if $\cT$ is infinite.
\end{thrm}
\begin{proof} We have already observed that if $\cT$ is finite,
then the trajectory must have a limit and this limit is a critical point of $G$
by Proposition \ref{pr.stationary}. Let us now show that the limit is a local
maximum of $G$. Arguing by contradiction, let us assume that $\cT$ is finite, and that the
limit is a local minimum of $G$.
Without loss of generality, we may assume that
$\lim_{t \to \infty} x(t) = 0$ and that $x$ is non-increasing for all sufficiently large times. Let $ \eps
>0$ be such that $g(\xi) \geq 2\,\delta \xi$ for $\xi\in ]0,\eps[$ and let $T>0$ be such that
$0 < x(t) < \eps$ for $t > T$. As the map $a$ converges to 0, one can choose $T$ to get
$a(t) < 2\,\sqrt{\delta}$ for $t>T$. Set $A(t) =
\exp\left(\frac{1}{2}\int_0^ta(s) ds\right)$ and $z(t) = A(t)x(t)$, then
\[ \ddot z(t) + \tilde g(t,z(t)) = 0 \]
where $$\tilde g(t,\xi) = A(t)g\left(\frac{\xi}{A(t)}\right) -
\frac{a(t)^2}{4}\xi - \frac{\dot a(t)}{2} \xi.$$
Recalling that $\dot a (t)\leq 0$ for every $t\geq 0$, we obtain
$\tilde g(t,\xi) \ge  \delta  \xi$, for $0 \le \xi \le A(t)\eps$. We derive that
$\ddot z(t) = - \tilde
g(t,z(t))\le 0$ and $z$ must be concave down. But since the function $z(t) = A(t) x(t)$ is also
positive for $t>T$, it must be increasing for all $t> T$, implying
 $\ddot z(t) \le -\delta z(t) \le -\delta z(T)$ for all $t >T$. This contradicts the fact that $z$ remains positive.

\medskip
We next consider the case where $\cT$ is infinite.
Without loss of generality, we can assume that $x(t_1) < x(t_2)$, and since $\dot{x}$ changes its
sign at each $t_k$, one sees that
\[x(t_{2j-1}) < x(t_{2j}), \quad x(t_{2j+1}) < x(t_{2j}) \]
for all $j \ge 1$. In fact, for all $j$ we have
\[x(t_{2j-1}) < x(t_{2j+1})  <  x(t_{2j+2}) < x(t_{2j})  \, .\]
Indeed, if \textit{e.g.} $x(t_{2j+2}) \ge x(t_{2j})$ for some $j$, then there
exists $s \in (t_{2j+1},t_{2j+2})$ (consequently $s > t_{2j}$) satisfying $x(s)
= x(t_{2j})$. And thus, since $\dot{x}(t_{2j})=0$, we have $w(s) \ge G(x(s)) =
G(x(t_{2j})) = w(t_{2j})$, which is a contradiction to the fact that $w$ is
strictly decreasing.

\medskip
Consequently, $X_1 = \lim_{j \to \infty} x(t_{2j+1})$ and $X_2 =
\lim_{j \to \infty} x(t_{2j})$ both exist, and $X_1 \le X_2$. We
claim that $X_1 = X_2$, which will prove that $\lim_{t \to \infty}
x(t) = X_1 = X_2$ exists. This limit must be a critical point of $G$, by
Proposition \ref{pr.stationary}. Since we have found a sequence $(x(t_k))_{k\geq 0}$
converging to it with $G(x(t_k)) > G(X_1) = G(X_2)$, the critical point can
only be a local minimum, completing the proof of the theorem.

\medskip
Suppose therefore $X_1 < X_2$. Clearly, $\lim_{t \to \infty}w(t)=G(X_1) = G(X_2)$
 since we have $\dot x(t_{2i}) =\dot x(t_{2i+1}) = 0$. From Proposition \ref{pr.accumulation},
 there exists a critical point $x^*$ of $G$ such that $\lim_{t \to \infty} w(t) = G(x^*)$. Since the trajectory does not converge, we
 deduce from Proposition \ref{pr.accumulation} that $x^*$ is not a local minimum of $G$.
Thus in view of assumption (c),  $x^*$ is a local maximum of $G$. Since $x^*$ is an accumulation point of the trajectory $x(.)$, we have $x^*\in [X_1, X_2]$. Observing that the sequence
$(x(t_{2j+1}))_{j\geq 0}$ converges to $X_1$ with $G(x(t_{2j+1}))>G(X_1)$, the point $X_1$ cannot
be a local maximum of $G$, hence $x^*\neq X_1$. The same argument shows that $x^*\neq X_2$,
and finally $x^*\in (X_1, X_2)$. Since $G(X_1)=G(X_2)=G(x^*)$, $X_1$ and $X_2$
 cannot be critical points of $G$ in view of assumption (b). Using the fact that $w$ is
 non-increasing, we have $G(\xi)\leq G(X_1) = G(X_2)$ for every $\xi \in [X_1, X_2]$.
We then deduce that  $G'(X_1)< 0$ and $G'(X_2)> 0$.  We are now in the situation of Lemma \ref{lemma_a},
with the limit points $X_1$ and $X_2$ coinciding respectively with the values defined by
(\ref{eq.def_xi1}) and (\ref{eq.def_xi2}). Hence if $t_i$ is sufficiently large, then
 by inequality (\ref{time})
\[t_{i+1} \le t_i + C + C \ln (t_i + 1) \]
with some constant $C$. By induction then
\[t_n \le C + C  n \ln (n+1) \]
for a suitable positive constant and for all sufficiently large $n$, say $n \ge N$.  These estimates imply
\begin{eqnarray*}
\int_{t_N}^\infty a(t)|\dot{x}(t)|^2 dt &\ge& \sum_{n \ge N}
\frac{c}{t_{n+1}+1} \int_{t_n}^{t_{n+1}}|\dot{x}(t)|^2dt \\
&\ge& \sum_{n \ge N} \frac{c D}{t_{n+1}+1} \\
&\ge& \sum_{n \ge N}
\frac{c D}{1+C+C\,(n+1) \ln (n+2)} = \infty \,
\end{eqnarray*}
where (\ref{dissipation}) in Lemma \ref{lemma_a} was used to estimate the
integrals $\int_{t_n}^{t_{n+1}}|\dot{x}(t)|^2dt$ from below by $D$. On the other hand,
$\int_{t_N}^\infty a(t)|\dot{x}(t)|^2 dt$ must be finite. This contradiction
proves the theorem.
\end{proof}

\begin{rmrk}[Assumptions on the map $a$]
A careful examination of the above proof shows that it is possible to slightly weaken the assumption
$a(t)\geq \frac{c}{t+1}$ for every $t\geq 0$. In fact, if we merely assume that
$$\int_1^\infty a(t\,\ln t) \, dt=\infty,$$
then the conclusions of Theorem \ref{thm_a} still hold true. We let the reader check that the map
$a$ defined by $t\mapsto \frac{1}{(t+1)\,\ln (\ln (t+3))}$ satisfies the above condition.
\end{rmrk}

\begin{rmrk}[Assumptions on the map $G$]
 Assumption (b) can be dropped, at the expense of a more technical proof.
 If assumption (c) is weakened, the current proof breaks down. For
example, if we merely assume that $G(x) \le G(x^*) - \delta\,|x-x^*|^p$ near local
maxima $x^*$ with some $p>2$, then we can only show that $t_{n+1} \le t_n + C +
C (1 + t_n)^{\frac{1}{2} - \frac{1}{p}}$, and stronger assumptions for $a$ must
be made, e.g. $a(t) \ge c (1+t)^{-\alpha}$ with $\alpha \left(\frac{3}{2} -
\frac{1}{p}\right) \le 1$. On the other hand, solutions may converge to local minima without
oscillating infinitely often, if $G$ is not strongly convex there. For example, $x(t) = (t+1)^{-\beta}$ is a solution of $\ddot{x}(t) + \frac{c}{t+1} \dot{x}(t) + x(t)^{1+ \frac{2}{\beta}} = 0$,
with $\beta>0$ and $c = 1 + \beta + \beta^{-1} $.
\end{rmrk}

\begin{rmrk}
Under the assumptions of Theorem \ref{thm_a}, for any solution of $(\cS)$ that converges to a local maximum $x^*$ of $G$, the set of sign changes $\cT = \{ t_1, \dots, t_K\}$ is finite. It appears plausible that $K$ and $t_K = \max \cT$ are bounded in terms of $\cE(0) = G(x(0)) + \frac12 |\dot{x}(0)|^2$, the potential $G$, and the function $a$.
\end{rmrk}

\medskip
We now show that under the assumptions of Theorem \ref{thm_a}, solutions generically converge to
a local minimum of $G$.

\begin{thrm} \label{thm_b}
 Under the assumptions of Theorem \ref{thm_a}, the set of initial data $(x_0,\,x_1)$
for which $\lim_{t \to \infty} x(t)$ is a local minimum of $G$ is open and dense.
\end{thrm}

\begin{proof}For $T>0$, define the map $F_T: \R^2 \to \R^2$ as $F_T(u,v) =
(x(T),\dot{x}(T))$, where $x$ is the solution of $(\cS)$ with $x(0)
= u, \, \dot{x}(0) = v$. By standard results for ordinary differential equations, see e.g.
\cite{Hirsch}, $F_T$ is a diffeomorphism and has an inverse $F_{-T}$. The inverse diffeomorphism
maps $(u,v) = (x(T),\dot{x}(T))$ to $F_{-T}(u,v) = (x(0),\dot{x}(0))$ by solving $(\cS)$
 backwards on $[0,T]$.

\medskip Let $x$ be a
solution for which $\lim_{t \to \infty} x(t) = \barre{x}$ is a local minimum
of $G$, with $x(0) = x_0$ and $\dot{x}(0) = x_1$. We shall find a
neighborhood of $(x_0, \, x_1)$ such that solutions with initial data from this
neighborhood have the same limit. There exist an open interval $I$ containing
$\barre{x}$ and $\delta > 0$ such that $\barre{x}$ is the only minimum of $G$ in $I$ and
$I$ is one of the connected components of $ G^{-1}\left([G(\barre{x}), G(\barre{x})+\delta)
\right)$. There is a time $T>0$ such that $x(t) \in I$ and $w(t) = G(x(t)) +
\frac{1}{2}|\dot{x}(t)|^2 < G(\barre{x})+\delta $ for $t>T$. Consider the open set
$\cO = \{(u,v) \, | u \in I,\, |v| < \sqrt{2 \delta + 2G(\barre{x}) - 2 G(u)} \, \}$.
By construction this set contains $(x(T),\dot{x}(T))$. Any solution $y$ of
$(\cS)$ with data $(y(T),\dot y(T)) \in \cO$ satisfies
$$
\forall t \ge T \qquad G(y(t)) \le w(t) \le w(T) \le G(\barre{x})+\delta.
$$
Using the definition of $\delta$ and $I$, we conclude that $y$ stays in $I$ for
all time greater than $T$. Since $I$ contains only one critical point of $G$
which is a local minimum, we infer that $y(t) \to \barre{x}$ as $t\to \infty$. Then $F_{-T}(\cO)$ is an
open neighborhood containing $(x_0, x_1)$, and all trajectories with initial
data in $F_{-T}(\cO)$ also converge to $\barre{x}$.

\medskip
Next let $\cI$ be the set of initial data in $\R^2$ whose solutions converge
to a local maximum of $G$. We must show that it has empty
interior. We first show the following:

\begin{clm} \label{cl.claimA}
Let $x$ be a solution with initial data $(x(0),\dot{x}(0)) \in \cI$ and let $x^*$ be
the limit of $x(t)$ as $t\to \infty$. In any neighborhood $\cU$ of $(x(0),\dot{x}(0))$,
 there exist $(y_0,y_1)$ such that the corresponding solution $y$ of ($\cS$) satisfies
$y^*=\lim_{t \to \infty} y(t)\neq x^*$ and $G(y^*)<G(x^*)$.
\end{clm}
 Let $\barre{\lambda} = \min \{ G(z)  |g(z) = 0, \, G(z) > G(x^*)\}$ be the next smallest critical value of $G$. Since $G$ is strictly concave near $x^*$, we can find
an interval $I = [x^*-\eps,\, x^* + \eps]$ such that
$g$ is strictly decreasing on $I$. Let $T\geq 0$ be such that $|x(T)- x^*| = \eps$ and $\dot{x}$ does
not change sign on $[T,\infty)$. After reducing $\eps$, we may assume that $w(T) < \barre{\lambda}$. After replacing $x$ with $x^*-x$
and $G(\xi)$ with $G(x^*-\xi)$ if necessary, we may also assume that
$x(T) = x^*-\eps < x(t) < x^*$ and $\dot{x}(t) > 0$ for $t > T$.
We claim that there exists a non-empty open interval $J$ containing $\dot{x}(T)$ such that whenever $y$ is a solution of $(\cS)$ with $y(T) = x(T)$ and $\dot{y}(T) \in J$, $\dot{y}(T) \ne \dot{x}(T)$,
then $y^* = \lim_{t \to \infty}y(t) \ne x^*$ and $G(y^*) < G(x^*)$.

\medskip
Indeed, first assume that there exists a solution $y\neq x$
 with $y(T) = x(T)$ that converges to $x^*$, such that $\dot{y}(t)>0$ for $t>T$. Then $v = x - y$ satisfies
\[\ddot v(t) + a(t)\dot v(t) + g(y(t)+v(t)) - g(y(t)) = 0 \]
as well as $v(T) = 0, \, \lim_{t \to \infty} v(t) = 0$. If $v$ has a positive maximum at some $t^*
> T$, then $\dot v(t^*) = 0$ and $\ddot v(t^*) \le 0$, hence $g(y(t^*) + v(t^*)) -
g(y(t^*)) = g(x(t^*)) - g(y(t^*)) \geq 0$. But we have $g(x(t^*)) <g(y(t^*))$ since
$x^*-\eps\leq y(t^*)<x(t^*)\leq x^*$ and since the map $g$ is decreasing on $[x^*-\eps,x^*]$,
a contradiction. The same argument applies if $v(t^{**})$ is a negative minimum of $v$.
So for any solution $y \ne x$ of $(\cS)$ with $y(T) = x(T)$ that converges to $x^*$,
 the derivative $\dot{y}$ must have at least one change of sign.

\medskip
Next let us assume that there is a sequence of solutions $y_k$ such that $y_k(T) = x(T)$,
 $\lim_{k\to\infty}\dot{y}_k(T)= \dot{x}(T)$, and $\lim_{t \to \infty} y_k(t) = x^*$ for all $k$. By the previous argument, the derivatives $\dot{y}_k$ must all change sign at least once on $(T,\infty)$. That is, for each $k$ there exists some minimal $t_k>T$ such $\dot{y}_k(t_k) = 0$. Then
$G(y_k(t_k)) > G(x^*)$ and hence $y_k(t_k) > x^* + \eps$. Let $T_k \in (T,t_k)$ be such that $y_k(T_k) = x^* + \eps$. By Remark \ref{re.majo}, especially inequality (\ref{time_est4}), we see that $T_k < t_k
\le T + C + C \ln (1+T)$ for all $k\in \N$, for some constant $C$.
By standard results on the continuous dependence of solutions of ordinary differential equations
 on initial data, we have
\begin{equation}\label{cont_dep}
\lim_{k \to \infty} \, \sup_{t \in [T,S]} |y_k(t) - x(t)| = 0
\end{equation}
for any $S\geq T$.  Recalling that $x(t)<x^*$ for every $t\geq T$, we have
$|y_k(T_k) - x(T_k)|>\eps$, which contradicts formula (\ref{cont_dep}) applied with
$S=T + C + C \ln (1+T)$. The contradiction shows that for some open interval $J$ containing $\dot{x}(T)$, solutions $y$ of $(\cS)$ with $y(T) = x(T)$ and $\dot{y}(T) \in J, \, \dot{y}(T) \ne \dot{x}(T)$ always have limits $y^*\ne x^*$. By shrinking the interval $J$ if necessary, we can guarantee that  $G(y^*)< \barre{\lambda}$ and therefore also $G(y^*)< G(x^*)$ for all such solutions. Consider then the set $F_{-T}\left(\{x(T)\} \times J\right)$. It contains
initial data arbitrarily close to $(x(0),\dot x(0))$ whose
solutions converge to a limit $y^*$ with $G(y^*) < G(x^*)$. This proves Claim \ref{cl.claimA}.

\medskip
To complete the proof of the theorem, let $n\geq 0$ be the number of local maxima of $G$.
If $G$ has no local maximum, the set $\cI$
is empty. Let us now assume that $n\geq 1$. Let $X_0=(x(0),\dot x (0))\in
\cI$ and let us denote by $x^*$ the limit of the corresponding solution $x$ of $(\cS)$. Let us fix
some $\varepsilon >0$. From Claim \ref{cl.claimA}, there exists $Y_0^{(1)}\in \R^2$ such that
$|Y_0^{(1)}-X_0|<\varepsilon/n$ and such that the limit $y^{(1),*}$ of the corresponding solution
of $(\cS)$ satisfies $G(y^{(1),*})<G(x^*)$. If $Y_0^{(1)}\in \cI$, we can apply again Claim \ref{cl.claimA}.
In fact, a repeated application of Claim \ref{cl.claimA} shows that there exist $k\leq n$ along with
$Y_0^{(2)}, \hdots, Y_0^{(k)}\in \R^2$ such that $Y_0^{(k)}\not \in \cI$ and
$$\forall i\in \{1,\hdots, k-1\}, \quad Y_0^{(i)} \in \cI \quad \mbox{ and } \quad
|Y_0^{(i+1)}-Y_0^{(i)}|\leq \varepsilon/n.$$
By summation, we derive that $|Y_0^{(k)}-X_0|\leq \frac{k\,\varepsilon}{n}\leq \varepsilon.$ Since the existence
of such a point $Y_0^{(k)}\not \in \cI$ is satisfied for every $\varepsilon >0$ and every $X_0\in
\cI$, we conclude that $\cI$ has empty interior.
\end{proof}

Let us finally establish a result that was used several times throughout this section.
Consider a coercive function $G:\R\to \R$ of class $\cC^1$ satisfying
assumptions (a)-(b) and let $x^*$ be a local maximum of $G$.
Let us define $X_1$, $X_2$ respectively by
\begin{equation}\label{eq.def_xi1}
X_1=\sup \{x\le x^*\,|\, G(x)>G(x^*)\},
\end{equation}
and
\begin{equation}\label{eq.def_xi2}
X_2=\inf \{x\ge x^*\,|\, G(x)>G(x^*)\}.
\end{equation}
The coercivity of $G$ shows that $-\infty<X_1\leq X_2<\infty$. Since $x^*$
is a local maximum, it is clear
that $X_1<x^*<X_2$. The continuity of $G$ shows that $G(X_1)=G(X_2)=G(x^*)$.
In view of assumption (b), this implies that $X_1$ and $X_2$ are not critical points of $G$.
Since $G(x)\leq G(X_1)=G(X_2)$ for every $x\in [X_1,X_2]$, we then have $G'(X_1)< 0$
and $G'(X_2)> 0$.

\begin{lmm} \label{lemma_a}
Let $a:\R_+ \to \R_+$ and $G:\R\to \R$ be as in Theorem \ref{thm_a}. Let $x^*$ be a
local maximum of $G$ and let $X_1$, $X_2$ be the real numbers respectively defined by
(\ref{eq.def_xi1}) and (\ref{eq.def_xi2}). Let $x(\cdot)$ be a
solution of $(\cS)$ and let $\cT = \{t_i\, |\, i \ge 1\}$ be the set of sign changes of $\dot{x}$.
Assume that for some $i \ge 1$,  $x(t_i)<X_1<X_2<x(t_{i+1})$ and $G'$ is negative
on $[x(t_i),X_1]$ and positive on $[X_2, x(t_{i+1})]$. Then there exist $C > 0, \, D>0$, and $T_0>0$ (defined in (\ref{T_0})) that depend only on $G, \, a,$ and the initial data such that if $t_i \ge T_0$, then
\begin{eqnarray}\label{dissipation}
\int_{t_i}^{t_{i+1}} |\dot{x}(t)|^2 dt &\ge& D \\
\label{time}
t_{i+1} - t_i &\le& C + C \ln (1+ t_i) \, .
\end{eqnarray}
These conclusions also hold true in the symmetric situation corresponding to
$x(t_{i+1})<X_1<X_2<x(t_{i})$.
\end{lmm}

\begin{proof} The assumption $x(t_i) < x(t_{i+1})$ implies
$\dot{x} > 0$ on $(t_i, \, t_{i+1})$. Since $G$ is strongly concave near $x^*$,
there exist $\eps$, $\delta >0$ such that $G(x) \le G(x^*) - \delta|x-x^*|^2$
whenever $|x-x^*|\le \eps$. We can assume that $x^* = 0$, $G(x^*) = 0$ and
that $X_1<-\eps<\eps<X_2$. Define numbers $\tau, s, \, s+h, \, \tau' \in (t_i, \, t_{i+1})$ by the conditions
$x(\tau) = X_1, \, x(s) = -\eps, \, x(s+h) = \eps, \, x(\tau') = X_2$.
We first prove inequality (\ref{dissipation}). Recalling that
 $w(t) \geq w(t_{i+1}) = G(x(t_{i+1}))>0$ for every  $t\in [t_i,t_{i+1}]$, we have
\begin{eqnarray*}
\int_{t_i}^{t_{i+1}}|\dot{x}(t)|^2 dt &=&\int_{t_i}^{t_{i+1}} \dot{x}(t)
\sqrt{2w(t) - 2G(x(t))} dt \\
&\ge& \int_{\tau}^{\tau'} \dot{x}(t) \sqrt{- 2G(x(t))} dt\\
&=& \int_{X_1}^{X_2} \sqrt{- 2G(x)} dx = D >0\, .
\end{eqnarray*}

\smallskip
We next prove inequality (\ref{time}). We first show that
\begin{equation}\label{time_est1}
s - t_i \le C, \,
\end{equation}
\begin{equation}\label{time_est2}
h \le C + C \ln (1+ t_{i+1}),
\end{equation}
and
\begin{equation}\label{time_est3}
t_{i+1} - (s+h) \le C\, .
\end{equation}
We claim that we can find $c>0$ such that for sufficiently large $t_i$,
\[w(t) \ge G(x(t_i)) - c(x(t)-x(t_i)) \ge G(x(t))\]
for $t \in [t_i,\tau]$. Indeed, choose  $0 < c < - \max_{[x(t_i),X_1]}G'$. Then $G(x)-G(x(t_i)) \leq
- c(x-x(t_i)) $ for $x \in (x(t_i),X_1)$. Now define
\begin{equation}\label{T_0}
T_0 = \inf \{s\geq 0 \, | a(s) \sqrt{2w(0)- 2 \min G} \le c\} \, .
\end{equation}
Since
$\dot{w}(t) = -a(t) |\dot{x}(t)|^2$ and $|\dot{x}(t)| \le \sqrt{2w(0) - 2\min G}$ for every
 $t\geq 0$, the inequality $\dot{w}(t) \ge -c |\dot{x}(t)|$ then holds for all
$t \ge T_0$. An integration then shows immediately that $w(t) \ge G(x(t_i)) - c(x(t)-x(t_i))$
 for $t \in [t_i,t_{i+1}]$ if $t_i \ge T_0$. For such $t_i$,
\begin{eqnarray*}
\tau - t_i &=& \int_{t_i}^{\tau} \frac{\dot{x}(t)}{\dot{x}(t)}dt =
\int_{t_i}^{\tau} \frac{\dot{x}(t)}{\sqrt{2w(t) - 2G(x(t))}} dt \\
&\le& \int_{t_i}^{\tau} \frac{\dot{x}(t)}{\sqrt{2G(x(t_i)) - 2c(x(t)-x(t_i)) - 2G(x(t))}}dt \\
&=& \int_{x(t_i)}^{X_1} \frac{dx}{\sqrt{2G(x(t_i)) - 2c(x-x(t_i)) - 2G(x)}}.
\end{eqnarray*}
The term under the square root is equivalent to $-2(G'(x(t_i))+c)\, (x-x(t_i))$ as $x\to x(t_i)$
and it ensues that the above integral is convergent, due to the choice of $c$.
Therefore, we derive
that $\tau - t_i \leq C_1$ for some constant $C_1$ that depends on $G$ and $c$. We next estimate $s - \tau$. Note that by construction, $G(x) < 0$ on $(X_1, -\eps]$ and $G'(X_1) < 0$. Hence $x \mapsto \frac{1}{\sqrt{-G(x)}}$ is integrable on $(X_1,-\eps]$. Then as before
\begin{eqnarray*}
s - \tau &=& \int_{\tau}^s \frac{\dot{x}(t)}{\sqrt{2w(t) - 2G(x(t))}} dt \\
&\le& \int_{\tau}^s \frac{\dot{x}(t)}{\sqrt{ - 2G(x(t))}}dt \\
&=& \int_{X_1}^{-\eps} \frac{dx}{\sqrt{- 2G(x)}} \\
&=& C_2
\end{eqnarray*}
for another constant $C_2$. These two estimates imply (\ref{time_est1}). Note that these two estimates do not use any information about $G$ on $[-\eps,x(t_{i+1})]$, and they hold also if in fact $x(\cdot)$ is monotone on $[t_i, \infty)$ and converges to $x^*$. That is,
\[s = \inf \{ t > t_i | x(t) \ge -\eps \} \le t_i + C\]
for some $C$ that depends only on $G$ and $a$.

Next, let us show (\ref{time_est3}), employing the same argument that was used to show (\ref{time_est1}):
\begin{eqnarray*}
t_{i+1} - (s+h) &=& \int_{s+h}^{t_{i+1}} \frac{\dot{x}(t)}{\sqrt{2w(t) - 2G(x(t))}} dt \\
&\le& \int_{s+h}^{t_{i+1}} \frac{\dot{x}(t)}{\sqrt{2G(x(t_{i+1})) - 2G(x(t))}}dt \\
&=& \int_{\eps}^{x(t_{i+1})} \frac{dx}{\sqrt{2G(x(t_{i+1}))- 2G(x)}} \\
&=& C_3
\end{eqnarray*}
where $C_3$ depends on $G$ and $\eps$.

We now show (\ref{time_est2}). Recall that $G(x) \le - \delta |x|^2$ on the interval $[-\eps,\eps]$ and that $w(t) \ge w(s+h)$ for every $t\in [s,s+h]$.
 Then
\begin{eqnarray*}
h &=& \int_{s}^{s+h} \frac{\dot{x}(t)}{\sqrt{2w(t) - 2G(x(t))}} dt \\
&\le& \int_{s}^{s+h} \frac{\dot{x}(t)}{\sqrt{2w(s+h) + 2\delta
|x(t)|^2}} dt \\
&=& \int_{-\eps}^{\eps}
\frac{dx}{\sqrt{2w(s+h) + 2\delta x^2}}\\
&=& \sqrt{\frac{2}{\delta}}\, \left[ \ln\left(x+\sqrt{x^2+\frac{w(s+h)}{\delta}} \right)\right]_0^\eps \\
&\le& C_4 -C_5\, \ln w(s+h),
\end{eqnarray*}
for suitable constants $C_4$, $C_5$ that depend on $\eps$ and $\delta$.
Let us estimate the quantity $w(s+h)$:
\begin{eqnarray*}
w(s + h) &=& w(t_{i+1}) + \int_{s+h}^{t_{i+1}}
a(t)|\dot{x}(t)|^2 dt \\
&\ge&\int_{s+h}^{t_{i+1}} \frac{c}{t+1}|\dot{x}(t)|^2 dt \ge
\frac{c}{t_{i+1}+1} \int_{s+h}^{t_{i+1}}|\dot{x}(t)|^2 dt.
\end{eqnarray*}
By arguing as in the proof of inequality (\ref{dissipation}), it is immediate to check that
$$\int_{s+h}^{t_{i+1}}|\dot{x}(t)|^2 dt\geq \int_{\eps}^{X_2} \sqrt{- 2G(x)} dx=D'>0.$$
It ensues that $ w(s+h)\geq \frac{c D'}{t_{i+1}+1}$ and consequently
\[h \le C + C \ln ( t_{i+1}+1) \]
for some $C>0$. Combining (\ref{time_est1}), (\ref{time_est2}) and (\ref{time_est3}) results in
\[ t_{i+1} - t_i \le C + C \ln (1+t_{i+1})\]
and therefore by an elementary argument
\begin{equation*}
t_{i+1} - t_i \le C + C \ln (1+t_i)
\end{equation*}
with some new constant $C>0$.  This proves (\ref{time}) completely.
\end{proof}

\begin{rmrk}\label{re.majo}
By combining (\ref{time_est2}) and (\ref{time_est3}), we obtain
\[t_{i+1} - s \le C + C \ln (1+t_{i+1})\]
and therefore by the same argument as above
\begin{equation}\label{time_est4}
t_{i+1} - s \le C + C \ln(1+ s),
\end{equation}
with some new constant $C>0$. Note that the proof of (\ref{time_est4}) does not use any
assumptions about the behavior of $x(\cdot)$ for $t<s$, which allows us to use it in the
proof of Theorem \ref{thm_b}.
\end{rmrk}

\section*{Appendix A: Stochastic approximation of $(\cS)$} \label{app.A}
We show in this appendix how the stochastic approximation scheme defined in the
introduction by (\ref{sa-modified}) naturally yields an ordinary differential equation
 of type $(\cS)$.
Using the same notations as in the introductory paragraph dealing with the
stochastic approximation, we define $h^{n+1}$ as the average drift at step $n$
by
$$
h^{n+1} =\frac{\displaystyle\sum_{i=0}^{n}
        \eps_i g(X^i,\omega^{i+1})}{\displaystyle\sum_{i=0}^{n} \eps_i}
$$
and we set $\tau^n=\eps_0+\eps_1 + \dots + \eps_n$. Note that $h^1$ is thus
initialized as $g(X^0,\omega^1)$. We can then rewrite (\ref{sa-modified}) in
\begin{equation}\label{rec}
\left\{
\begin{array}{ll}
        (X^0,\tau^0,h^0) \in \R^d \times \{\eps_0\}\times \{0\}\\[0.5cm]
        \tau^{n+1} = \displaystyle\tau^n + \eps_{n+1}\qquad \forall n \in \N\\[0.5cm]
        h^{n+1} = \displaystyle h^n - \eps_{n} \frac{h^n}{\tau^{n}}+\eps_{n}
        \frac{g(X^{n},\omega^{n+1})}{\tau^{n}}\qquad \forall n \in \N\\[0.5cm]
        X^{n+1} = \displaystyle X^n-\eps_{n+1} h^{n+1}\qquad \forall n \in \N.
    \end{array}
\right.
\end{equation}
The recursion (\ref{rec}) can be identified now as
$$(\tau^{n+1},h^{n+1},X^{n+1}) = (\tau^{n},h^{n},X^{n}) - \eps_{n+1}
H(\tau^n,h^n,X^{n},\omega^{n+1}) + \eps_{n+1} \eta^{n+1}$$ where
$H(\tau,h,X,\omega)=(-1,(h-g(X,\omega))/\tau,h)$
 and the residual perturbation $\eta^{n+1}$ is
given by
$$
\eta^{n+1} = \left(0,(\eps_{n+1}-\eps_n)
\frac{h^n-g(X^n,\omega^{n+1})}{\tau^n\eps_{n+1}},\eps_{n}\frac{h^n-g(X^n,\omega^{n+1})}{\tau^{n}}\right).
$$
Remark then that $\eta^{n} \to 0$ as $n \to \infty$ and the
conditional expectation with respect to filtration $\mathcal{F}_n$ of
$H(\tau^n,h^n,X^{n},\omega^{n+1})$ is
$$
\mathbb E \left[ H(\tau^n,h^n,X^{n},\omega^{n+1})| \mathcal{F}_n \right] =
\left( -1,\frac{h^n-g(X^n)}{\tau^n},h^n \right).
$$
Applying the result of \cite{robbins-monro}, the time interpolation of the
process $(\tau^n,h^n,X^n)_{n \geq 1}$ asymptotically behaves as the solution of the following
system of differential equations
\begin{equation}\label{sa-ode}
\left\{
\begin{array}{ll}
\dot{\tau}(t) = 1 \qquad \forall t \in \R\\[0.5cm]
\dot{h}(t) = -\frac{h(t)-g(X(t))}{\tau(t)} \qquad \forall t \in \R\\[0.5cm]
\dot{X}(t) = - h(t) \qquad \forall t \in \R.
    \end{array}
\right.
\end{equation}
If we note $\beta=\tau(0)$, we now have
$$
\ddot{X}(t)= - \dot{h}(t) = \frac{h(t)-g(X(t))}{t+\beta}  = -\frac{\dot{X}(t)+g(X(t))}{t+\beta},
$$
which is a particular case of $(\cS)$.

\section*{Appendix B: Special Cases} \label{app.B}
Consider first the equation
\begin{equation} \label{bessel_c}
\ddot{x}(t) + \frac{c}{t} \dot{x}(t) + x(t) = 0
\end{equation}
for $t > 0$ and its shifted version
\begin{equation} \label{bessel_c_1}
\ddot{x}(t) + \frac{c}{t+1} \dot{x}(t) + x(t) = 0 \, .
\end{equation}
Let $\tilde{x}$ be a solution of the Bessel equation
\[t^2\ddot{x}(t) + t \dot{x}(t) + \left(t^2 - \left(\frac{c-1}{2}\right)^2\right)x(t) = 0
\]
for $t > 0$, i.e.
\[\tilde{x}(t) = b_1J_{(c-1)/2}(t) +  b_2 Y_{(c-1)/2}(t)\]
where $J_{(c-1)/2}$ and $Y_{(c-1)/2}$ are Bessel functions of the
first and second kind. A calculation shows that
\[x(t) = t^{\frac{1-c}{2}} \tilde {x}(t) \]
is a solution of (\ref{bessel_c}). For $x(0)$ to be finite, we require $b_2 = 0$.
Hence the general solution of (\ref{bessel_c}) with finite $x(0)$ is
\[x_c(t) = b_1 t^{\frac{1-c}{2}}J_{(c-1)/2}(t) \, .\]
Since $J_{(c-1)/2}(t) \approx \sqrt{\frac{2}{\pi t}} \cos (t - \frac{c\pi}{4} )$ as $t \to \infty$ for all $c$, we therefore see that
\[x_c(t) \approx C t^{-c/2} \cos (t - \frac{c\pi}{4})\]
as $t \to \infty$, for all $t$, with a suitable constant C.

\medskip
Solutions of (\ref{bessel_c_1}) are of the form
\[x(t) = b_1 (t+1)^{\frac{1-c}{2}}J_{(c-1)/2}(t+1) + b_2 (t+1)^{\frac{1-c}{2}}Y_{(c-1)/2}(t+1) \]
and have the asymptotic behavior
\[x(t) \approx  t^{-c/2} \left( C \cos (t - \varphi_0) \right)
\]
with a suitable amplitude constant $C$ and phase shift $\varphi_0$.
This is the typical behavior of solutions of ($\cS$) in one
dimension near non-degenerate local minima of $G$, in the case $a(t)
= ct^{-1}$ or $a(t) = c(t+1)^{-1}$.

\medskip
Consider now the equation
\begin{equation} \label{bessel_c_mod}
\ddot{y}(t) + \frac{c}{t} \dot{y}(t) - y(t) = 0
\end{equation}
for $t > 0$ and its shifted version
\begin{equation} \label{bessel_c_mod1}
\ddot{y}(t) + \frac{c}{t+1} \dot{y}(t) - y(t) = 0 \, .
\end{equation}
We are interested in solutions $y_c$ that converge to $0$.
Let $\tilde{y}$ be a solution of the modified Bessel equation
\[t^2\ddot{y}(t) + t \dot{y}(t) - \left(t^2 + \left(\frac{c-1}{2}\right)^2\right)y(t) = 0
\]
for $t > 0$, i.e.
\[\tilde{y}(t) = b_1I_{(c-1)/2}(t) +  b_2 K_{(c-1)/2}(t)\]
where $I_{(c-1)/2}$ and $K_{(c-1)/2}$ are modified Bessel functions of the first and second kind.
A direct calculation shows again that
\[y_c(t) = t^{\frac{1-c}{2}} \tilde {y}(t) \]
is a solution of (\ref{bessel_c_mod}). For $y_c$ to be convergent to $0$, we require $b_1 = 0$,
since $I_\nu(t) \to \infty$ as $t \to \infty$. Thus solutions of (\ref{bessel_c_mod}) that converge to zero are of the form
\[y_c(t) =b_2 t^{\frac{1-c}{2}}K_{(c-1)/2}(t) \, .\]
Since $K_{(c-1)/2}(t) \approx \sqrt{\frac{2}{\pi t}}e^{-t}$ for all $c$ as $t \to \infty$
with higher order terms depending on $c$, we see that
\[y_c(t) \approx C t^{-c/2} e^{-t}\]
as $t \to \infty$, with some constant $C$. Solutions of (\ref{bessel_c_mod1}) that converge to $0$ then are of the form
\[y(t) = b_2 (t+1)^{\frac{1-c}{2}}K_{(c-1)/2}(t+1)\]
and have the same asymptotic behavior. This is the typical behavior
of solutions of ($\cS$) in one dimension near non-degenerate local
maxima of $G$, again in the case $a(t) = ct^{-1}$ or $a(t) =
c(t+1)^{-1}$. The standard reference for results on Bessel functions
is \cite{abramowitz}.

\end{document}